\documentclass[reqno]{amsart}

\usepackage{myPackages}
\usepackage{myCommands}
\usepackage{zetaCommands}

\usetikzlibrary{intersections, calc, positioning, arrows.meta}

\newcommand{\Li}{\mathrm{Li}}

\newcommand{\cH}{\mathcal{H}}
\newcommand{\cA}{\mathcal{A}}
\newcommand{\cU}{\mathcal{U}}
\newcommand{\cV}{\mathcal{V}}

\newcommand{\bk}{\mathbf{k}}

\newcommand{\tI}{\widetilde{\mathrm{I}}}

\begin{document}

\title{
Motivic interpretations for iterated integrals on some specific algebraic curves
}
\author{
Eisuke Otsuka
}
\date{}

\address{Eisuke Otsuka\\
Mathematical Institute, Tohoku University.\\
6-3, Aoba, Aramaki, Aoba-ku, Sendai, 980-8578, JAPAN}
\email{eisuke.otsuka.p3@dc.tohoku.ac.jp}

\maketitle

\begin{abstract}

  Multiple zeta values (MZVs for short) can be represented as iterated integrals of $\bQ$-rational algebraic differential forms on $\bP^1(\bC)\setminus\{0, 1, \infty\}$. This interpretation allows us to consider MZVs geometrically, and this is one of the motivations for Deligne--Goncharov, Terasoma et al. to give motivic interpretations of MZVs by using the theory of mixed Tate motives and the motivic fundamental groups.

  In this paper, we consider the iterated integrals on some rational curves over $\bQ$ and study their arithmetic properties. They are an extension of MZVs and also include some other known special values such as multiple $\widetilde{T}$-values. Furthermore, we give motivic interpretations of them by investigating a relationship with motivic iterated integrals given by Goncharov. At this point, it is important to consider the base expansion and the Galois invariant part of the space of motivic iterated integrals. Finally, we denote that a motivic interpretation of the alternating multiple mixed values can be given by the same method. Our results also extend a part of author's previous work.
\end{abstract}

\tableofcontents

\section{Introduction} \label{sec: introduction}
The multiple zeta values (MZVs, also called the Euler-Zagier type MZVs) are defined by 
\begin{align}
\zeta(k_1,\dots,k_d)=\sum_{0<n_1<\dots<n_d}\frac{1}{n_1^{k_1}\dots n_d^{k_d}}
\end{align}
for $k_1,\dots,k_{d-1}\in\bZge{1}$, and $k_d\in\bZge{2}$. A pair $\bk=(k_1,\dots,k_d)\in(\bZge{1})^d$ with $k_d>1$ is called the (admissible) index of the weight $|\bk|:=k_1+\dots+k_d$. It is known that there are many $\bQ$-linear relations along MZVs (\cite{IhKaZa06},\cite{LeMu95},\cite{Oh99} et al.), and in \cite{Za94}, Zagier conjectured that the dimension of the $\bQ$-linear space $Z_k$ spanned by MZVs of weight $k$ is
\begin{align}
    \dim_\bQ Z_k=d_k,
\end{align}
where
\begin{align} \label{eq: dimMZV}
\sum_{k\ge0}d_kt^k=\frac{1}{1-t^2-t^3}.
\end{align}
This conjecture has not yet been proved as of August 2024, but Goncharov (\cite{Go01}), Terasoma (\cite{Te02}), Deligne--Goncharov (\cite{DeGo05}) showed 

\begin{align}
    \dim_\bQ Z_k\le d_k
\end{align}
by using the category of the mixed Tate motives over $\bQ$. At that time, the iterated integral representation given below was important to give a motivic interpretation of MZVs. The iterated integral representation of MZVs is said to have been first mentioned by M. Kontsevich (see \cite[p.510]{Za94}) and is given by
\begin{align}
\zeta(k_1,\dots,k_d)=\int_{0<x_1<\dots<x_k<1}\prod_{j=1}^k\eta_j(x_j), \label{eq: EZint}
\end{align}
where $$\eta_j(x)=\begin{cases}
\omega_1=\omega_1(x):=\frac{dx}{1-x}, & j=k_1+\dots+k_{s-1}+1 \text{ for some } s=1,\dots,d,\\
\omega_0=\omega_0(x):=\frac{dx}{x}, & \text{otherwise}.
\end{cases}$$

We study the special values defined by iterated integrals. In general, for 1-forms $\eta_1,\dots,\eta_k$ on a manifold $M$ and a smooth path $\gamma: [0,1]\rightarrow M$, we define
\begin{align}
  \int_\gamma \eta_1\cdots\eta_k:=\int_{0<t_1<\cdots<t_k<1}\gamma^*\eta_1(t_1)\cdots\gamma^*\eta_k(t_k),
\end{align}
where $\gamma^*\eta_j~(j=1,\dots,k)$ is the pullback of $\eta_j$ along $\gamma$. For example, if $M=\bP^1(\bC)\setminus\{0,1,\infty\}$, we have the iterated integral representations of the multiple zeta values:
\begin{align} \label{eq: MZV}
  \zeta(k_1,\dots,k_r)=\int_\dch\eta_1\cdots\eta_k,
\end{align}
where $k=k_1+\cdots+k_r$, $\dch:[0,1]\rightarrow\bP^1(\bC),\dch(t):=t$, and 
$$
  \eta_j:=\begin{cases}
  \omega_1, & j=k_1+\cdots k_s+1 \text{ for some } s=0,\dots,r-1, \\
  \omega_0, & \text{otherwise}.
\end{cases}
$$

In this paper, we consider iterated integrals on some specific algebraic curves and give their motivic interpretations. We also study explicit formulas for some iterated integrals. This paper is also regarded as a sequel of our previous work (\cite{Ot23}).

\subsection{Main Results} \label{subsec: results}

Let $f\in\bQ[X,Y,Z]$ be an irreducible homogeneous polynomial such that $$X_f:=\left\{[X:Y:Z]\in\bP^2(\bC)~\middle|~f(X,Y,Z)=0\right\}$$ is a geometrically connected and projective smooth curve over $\bQ$.
We define $\varphi_f:X_f\rightarrow\bP^1(\bC), \varphi_f([X:Y:Z]):=[X:Z]$ and consider a local coordinates $x=\frac{X}{Z}$ and $y=\frac{Y}{Z}$ for $Z\ne0$. Put $Y_f:=X_f\setminus\{x=0,1,\infty\}$.
Let $\dch_f: [0,1]\rightarrow X_f$ be a smooth path which satisfies $\varphi_f\circ\dch_f=\dch$.
Fix a $\bQ$-basis $B_f$ of the first algebraic de Rham cohomology $H^1_{\mathrm{dR}}(Y_f/\bQ)$ (cf. \cite{Ke}).

\defi \label{def: II on Xf}
For $\eta_1,\dots,\eta_k\in B_f$, we define
\begin{align} \label{eq: def of II on Xf}
  \mathrm{I}_f(\eta_1\cdots\eta_k):=\int_{\dch_f}\eta_1\cdots\eta_k.
\end{align}
We consider only the case when the \cref{eq: def of II on Xf} converges and set 
\begin{align} \label{eq: Bfk}
  B_f^{(k)}:=\left\{(\eta_1,\dots,\eta_k)\in (B_f)^k~\middle|~\mathrm{I}_f(\eta_1\cdots\eta_k) \text{ converges}\right\}.
\end{align}
We call $k$ the weight of $(\eta_1,\dots,\eta_k)\in B_{f}^{(k)}$, and let $Z_{f}^{(k)}$ be the $\bQ$-linear space spanned by $\mathrm{I}_f(\eta_1\cdots\eta_k)$ for $(\eta_1,\dots,\eta_k)\in B_f^{(k)}$. When $k=0$, we interpret $\mathrm{I}_f(\eta_1\cdots\eta_k)=1$ and set $Z_f^{(k)}=\bQ$.
\edefi

In this paper, we consider the following case when $f\in\{g,h\}$ where
\begin{enumerate}[(I)]
  \item $g(X,Y,Z):=X^2+Y^2-Z^2;$
  \item $h(X,Y,Z):=X^2+XY+Y^2-Z^2.$
\end{enumerate}

Note that $X_f$ is defined over $\bQ$ with good reduction outside $S_f$ where
\begin{align}
    S_f=\begin{cases}
        \{2\}, & f=g,\\
        \{2,3\}, & f=h.
    \end{cases}
\end{align}

The main result of this paper is to give a motivic interpretation for these values by using motivic iterated integrals. To explain the claim, we need more notation. Fix an algebraic closure $\overline\bQ$ of $\bQ$. For a number field $F(\subseteq\overline\bQ)$ and the $\bQ$-vector subspace $R\subseteq F^\times\otimes_\bZ\bQ$ with $R=\mathcal{O}^\times_{F,S}\otimes_\bZ\bQ$ for some finite set $S$ of finite places of $F$, let $\cA_{F,R}$ be the fundamental Hopf algebra of $R$ and $\cH_{F,R}:=\cA_{F,R}\otimes_\bQ\bQ[\tau]$. Let $\per: \cH_{F,R}\rightarrow \bC$ be the period map. Then the Galois group $\Gal(F/\bQ)$ acts on $\cH_{F,R}\otimes_\bQ F$ (see \cref{sec: motivic} for details).

We apply them to our case for $f\in\{g,h\}$. We define a number field $F_f$ by
$$F_f:=\begin{cases}
  \bQ(\xi_4), & f=g,\\
  \bQ(\xi_6), & f=h,\\
\end{cases}$$
where $\xi_N:=\exp\left(\frac{2\pi\sqrt{-1}}{N}\right)$. We define
\begin{align}
    R_f:=\mathcal{O}_{F_f,S'_f}^\times\otimes_\bZ\bQ=\begin{cases}
        \langle2\rangle_{F_g}, & f=g,\\
        \langle2,3\rangle_{F_h}, & f=h,
    \end{cases}
\end{align}
where 
\begin{align}
    S'_f&:=\{v \text{ finite primes of } F_f~|~v|p \text{ for any }p\in S_f\}\\
    &=\begin{cases}
            \{(1-\xi_4)\}, & f=g,\\
            \{(2),(1-\xi_3)\}, & f=h
        \end{cases}
\end{align}
and $\langle a_1,\dots,a_r\rangle_F$ is the $F$-linear space spanned by $a_1,\dots,a_r$. 
Also, we put $G_f:=\Gal(F_f/\bQ)$. Then, we have the following theorem. 

\thm \label{thm: 1} 
For each $f\in\{g,h\}$ and $(\eta_1,\dots,\eta_k)\in B_f^{(k)}$, there exists an element $$\Im_f(\eta_1\cdots\eta_k)\in\left(\cH_{F_f,R_f}^{(k)}\otimes_\bQ F_f\right)^{G_f}$$
such that $\per(\Im_f(\eta_1\cdots\eta_k))=\mathrm{I}_f(\eta_1\cdots\eta_k)$.
\ethm

As an application of \cref{thm: 1}, we have an upper bounds of the dimension of $Z_f^{(k)}$ for each $f\in\{g,h\}$.

\cor \label{cor: dimension}
For $f\in\{g,h\}$, it holds
\begin{align}
  \dim_\bQ Z_{f}^{(k)}\le D_f^{(k)}
\end{align}
where the integer $D_f^{(k)}$ is defined by
\begin{align}
  \sum_{k=0}^\infty D_f^{(k)}t^k=\begin{cases}
    \frac{1}{1-2t}, & f=g,\\
    \frac{1}{1-3t+t^2}, & f=h.
  \end{cases}
  \end{align}
\ecor

We will compute the value of $D_f^{(k)}$ and conjectural value of $\dim_\bQ Z_f^{(k)}$ after \cref{cor: dim2} obtained by numerical experiments. 

Our (motivic) periods are coming from the motives defined over $\bQ$. However, if we extend the base field $\bQ$ suitably, they are regarded as the (motivic) periods coming from the category of the mixed Tate motives over cyclotomic fields which are well-known objects (\cite{Go05}). 
Thus, a key, as in previous paper \cite{Ot23}, is to consider the Galois descent of the (motivic) periods after a suitable base extension. 

We organize this paper as follows. In \cref{sec: preparations}, we set up some notations and definitions. In \cref{sec: basic}, we explain some basic properties of our iterated integrals and give some explicit formulas. These formulas shows our periods can be written as linear combinations of well known periods and it leads us to seek their motivic interpretations, which are the contents in \cref{sec: motivic,sec: results}. In \cref{sec: motivic}, we summarize the basic properties of the motivic iterated integrals given by Deligne, Glanois, and Goncharov to give motivic interpretations of the iterated integrals on our algebraic curves. In \cref{sec: results}, we give motivic interpretations and investigate their algebraic structures. This interpretation explains the appearance of $\pi$, $\log2$, $\log3$, $L(k,\chi_{-3})$, and Riemann zeta values in the explicit formulas in \cref{ex: example}-(3), \cref{prop: special-case1,prop: special-case2}. We also notice about the direct sum decomposition of the $\bQ$-linear space $Z_f^{(k)}$. In \cref{ss: AMMV}, we describe the relationship between our periods and special values called alternating multiple mixed values, defined by Xu--Yan--Zhao \cite{XuYaZh23} and further studied by Charlton \cite{Ch24}.

We remark that we considered in \cite{Ot23} only the values given by iterated integrals with three elements of $B_g$ and we consider a wider class in this paper. 

\section*{Acknowledgement}
The author would like to express his deepest gratitude to his supervisor, Professor T. Yamauchi. In particular, we follow his idea in a part of the computation in \cref{sec: basic} and Subsection \ref{subsec: Galoisact}. The author would also like to thank Professor H. Furusho of Nagoya University for inviting me to his monthly seminar and gave me various advice of my research. The author is grateful to Professor S. Charlton of the Max-Planck Institute for Mathematics for discussing alternating multiple mixed values. The author would like to thank the WISE Program for AI Electrics of Tohoku University for their financial support to carry out our study.

\section*{Notation}
Throughout this paper, we use the following notation unless otherwise mentioned.
\nota \label{notation}
~
\begin{itemize}
  \item For a finite set $S$, let $\#S$ be the cardinality of $S$.
  \item Kronecker delta: for $r,s\in\bZ$, put $\delta_{r,s}:=\begin{cases}
  1, & r=s,\\
  0, & r\ne s.
  \end{cases}$
  \item Binomial coefficients: for $k,r\in\bZge{0}$, put $\binom{k}{r}=\begin{cases}
  \frac{k!}{r!(k-r)!}, & k\ge r,\\
  0, & k<r.\end{cases}$
  \item For $n,r,s\in\bZge{1}$, let $S_n$ be the $n$-th symmetry group, and define $(r,s)$-shuffle $S_{r,s}$ by $$S_{r,s}:=\left\{\delta\in S_{r+s}~\middle|~\begin{array}{c} \delta^{-1}(1)<\dots<\delta^{-1}(r)\\ \delta^{-1}(r+1)<\dots<\delta^{-1}(r+s)\end{array}\right\}.$$
  \item For $N\in\bZge{1}$, put $\xi_N:=\exp\left(\frac{2\pi\sqrt{-1}}{N}\right)$.
  \item For $N\in\bZge{1}$, put $\mu_N:=\left\{\xi_N^r~\middle|~r=0,1,\dots,N-1\right\}$ and $\widetilde{\mu}_N:=\mu_N\cup\{0\}$.
  \item Path reversal: for a piecewise smooth path $\gamma:[0,1]\rightarrow\bC$, let $\gamma^{-1}:[0,1]\rightarrow\bC$ be $$\gamma^{-1}(t):=\gamma(1-t).$$
  \item Path connection: for piecewise smooth paths $\gamma_1, \gamma_2:[0,1]\rightarrow\bC$ satisfying $\gamma_1(1)=\gamma_2(0)$, let $\gamma_1\gamma_2: [0,1]\rightarrow\bC$ be $$(\gamma_1\gamma_2)(t):=\begin{cases}
  \gamma_1(2t), & 0<t<1/2,\\
  \gamma_2(2t-1), & 1/2<t<1.
  \end{cases}$$
  \item  For $p,q\in\bC$, we define the straight path $\dch_{p,q}$ on $\bC$ by $$\dch_{p,q}: [0,1]\rightarrow\bC,~~\dch_{p,q}(t)=p+t(q-p).$$ In particular, for $p=0$ and $q=1$, we put $\dch:=\dch_{0,1}$.
  \item  For $p,q_1,q_2\in\bC$ with $|p-q_1|=|p-q_2|=r$, we define the arc $C_p(q_1,q_2)$ on $\bC$ by $$C_p(q_1,q_2): [0,1]\rightarrow\bC,~~C_p(q_1,q_2)(t)=p+re^{\sqrt{-1}(\arg(q_1)+t(\arg(q_2/q_1)))}.$$
  \item For a number field $F$, we put $r_1(F):=\#\{\text{real places}\}$ and $r_2(F):=\#\{\text{complex places}\}$.
  \item For an $F$-linear space $V$, let $T(V)=\bigoplus_{k=0}^\infty V^{\otimes k}$ be the tensor algebra of $V$.
  \item For a field $F$ and a set $S$, let $F\langle S\rangle:=F\langle e_s~|~s\in S\rangle$ be the non-commutative $\bQ$-algebra generated by the formal basis $e_s~(s\in S)$.
  \item For a scheme $X$ with locally of finite type over $\bC$, let $X^\mathrm{an}$ be the analytification of $X$.
\end{itemize} 
\enota

\section{Preparations} \label{sec: preparations}
First, we prepare some notation for iterated integrals on $\bC$.

\defi \label{def: II on C}
Let $\gamma:[0,1]\rightarrow\bC$ be a piecewise smooth path on $\bC$ and $a_1,\dots,a_k\in\bC\setminus\gamma((0,1))$ with $p:=\gamma(0)\ne a_1, q:=\gamma(1)\ne a_k$. Then, we define
\begin{align} \label{eq: II on C}
\mathrm{I}_\gamma(p;a_1,\dots,a_k;q):=\int_\gamma\frac{\diff x}{x-a_1}\dots\frac{\diff x}{x-a_k}.
\end{align}
For the moment, we do not ask if the above iterated integral converges.

In addition, for each element $w=e_{a_1}\dots e_{a_k}\in \bQ\langle\bC\rangle$ with $a_1\ne p, a_k\ne q$, we define
\begin{align}
\mathrm{I}_\gamma(p;w;q):=\mathrm{I}_\gamma(p,a_1,\dots,a_k;q)
\end{align}
and extend $\mathrm{I}_\gamma(p;w;q)$ $\bQ$-linearly to any $w\in \bQ\langle\bC\rangle$.
\edefi

\rem \label{rem: reg of II}
The right hand side of the \cref{eq: II on C} diverges if $p=a_1$, $q=a_k$, or $a_1,\dots,a_k\in\gamma((0,1))$. However, we can regularize $\mathrm{I}_\gamma(p;a_1,\dots,a_k;q)\in\bC$ in these cases as in {\cite[p.6, Definition 2.0.6]{Ot23}}. For example, if $k=1$, then it holds
\begin{align}
  \exp(\mathrm{I}_{\dch_{a,c}}(a;b;c))=\tI(a,b,c),
\end{align}
where
\begin{align}
  \tI(a,b,c):=\begin{cases}
    \frac{c-b}{a-b}, & a\ne b \text{ and } b\ne c,\\
    c-b, & a=b \text{ and } b\ne c,\\
    (a-b)^{-1}, & a\ne b \text{ and } b=c,\\
    1 & a=b=c.
  \end{cases}
\end{align}
\erem

\defi
For $k_1,\dots,k_d\in\bZge{1}$ and $\alpha_1,\dots,\alpha_d\in\mu_N$ with $(k_d,\alpha_d)\ne(1,1)$, 
\begin{align}
  L\left(\begin{matrix}
    k_1,\dots,k_d\\
    \alpha_1,\dots,\alpha_d
  \end{matrix}\right):=(-1)^d\mathrm{I}_\dch(0;\alpha_1^{-1},\{0\}^{k_1-1},\dots,\alpha_d^{-1},\{0\}^{k_d-1};1)
\end{align}
is called the (shuffle type) multiple $L$-value of level $N$.
\edefi

Next, we consider iterated integrals in \cref{def: II on Xf} for our algebraic curves.

\exm \label{ex: II on Xf}
If $f(X,Y,Z)=X^2+Y^2-Z^2$, we can take a $\bQ$-basis as $B_f=\left\{\frac{\diff x}{x},\frac{\diff x}{1-x}, \frac{\diff x}{y}, \frac{\diff x}{xy}\right\}$. Then, if we restrict the $1$-forms to the subset $\left\{\frac{\diff x}{x},\frac{\diff x}{1-x}, \frac{\diff x}{y}\right\}\subseteq B_f$, the iterated integrals in the \cref{eq: def of II on Xf} are $F_2$-MZVs considered in \cite{Ot23}.
\eexm

Now, we consider the cases of $g(X,Y,Z)=X^2+Y^2-Z^2$ and $h(X,Y,Z)=X^2+XY+Y^2-Z^2$. Note that for each $f\in\{g,h\}$, the geometric genus of $X_f$ is $0$. Further, $X_g$ is smooth over $\bZ[1/2]$ and $X_h$ is smooth over $\bZ[1/3]$.

\vspace{5mm}
\noindent
{\textbf{An explicit $\bQ$-basis of $H_{\mathrm{dR}}^1(Y_f/\bQ)$.}}
For each $X_f$ ($f\in\{g,h\}$), we specify $B_f$ as follows.
\begin{enumerate}[(I)]
  \item \underline{$g(X,Y,Z)=X^2+Y^2-Z^2$}: There is an isomorphism
  \begin{align}
    \phi_g: \bP^1(\bC)\overset{\sim}{\longrightarrow}(X_g)^\mathrm{an}, \phi_g(\lambda):=(x,y)=\left(\frac{2\lambda}{1+\lambda^2},\frac{1-\lambda^2}{1+\lambda^2}\right).
  \end{align}
  The points on $X_g$ with $x=0,1,\infty$ correspond to the following five points: $(x,y)=\left(0,\pm1\right), (1,0)$ and $[X:Y:Z]=\left[1:\pm\sqrt{-1}:0\right]$. As mentioned in \cite[p.10]{Ke}, it holds
  \begin{align}
      H^1_\dR(Y_g/\bQ)\otimes_\bQ\bC\cong H^1_\dR((Y_g)^{\mathrm{an}})\cong H^1_\dR(\bP^1(\bC)\setminus\{\text{five points}\}),
  \end{align}
  then the dimension of the $\bQ$-linear space $H^1_\dR(Y_g/\bQ)$ is $4$ and we can take $B_g$ as $$B_g=\left\{\omega_0:=\frac{\diff x}{x},\omega_1:=\frac{\diff x}{1-x}, \omega_2:=\frac{\diff x}{y}, \omega_3:=\frac{\diff x}{xy}\right\}.$$
  
  The set $B_g^{(k)}$ defined by the \cref{eq: Bfk} becomes
  $$B_g^{(k)}=\{(\eta_1,\dots,\eta_k)\in (B_g)^k~|~\eta_1\ne\omega_0,\omega_3 \text{ and } \eta_k\ne\omega_1\}.$$
  \item \underline{$h(X,Y,Z)=X^2+XY+Y^2-Z^2$}: There is an isomorphism
  \begin{align}
    \phi_h: \bP^1(\bC)\overset{\sim}{\longrightarrow}(X_{h})^\mathrm{an}, \phi_{h}(\lambda):=(x,y)=\left(\frac{\lambda(\lambda+2)}{1+\lambda+\lambda^2},\frac{1-\lambda^2}{1+\lambda+\lambda^2}\right).
  \end{align}
  The points on $X_h$ with $x=0,1,\infty$ correspond to the following six points: $(x,y)=(0,\pm1), (1,0),(1,-1)$ and $[X:Y:Z]=[1:\xi_3^{\pm1}:0]$. Since
  \begin{align}
      H^1_\dR(Y_h/\bQ)\otimes_\bQ\bC\cong H^1_\dR((Y_h)^{\mathrm{an}})\cong H^1_\dR(\bP^1(\bC)\setminus\{\text{six points}\}),
  \end{align}
  the dimension of the $\bQ$-linear space $H^1_\dR(Y_h/\bQ)$ is $5$ and we can take $B_h$ as $$B_h=\left\{\omega_0:=\frac{\diff x}{x},\omega_1:=\frac{\diff x}{1-x}, \omega_4:=\frac{\diff x}{x+2y}, \omega_5:=\frac{\diff x}{x(x+2y)},\omega_6:=\frac{\diff x}{(1-x)(x+2y)}\right\}.$$ The set $B_h^{(k)}$ defined by the \cref{eq: Bfk} becomes
  $$B_h^{(k)}=\{(\eta_1,\dots,\eta_k)\in (B_h)^k~|~\eta_1\ne\omega_0,\omega_5 \text{ and } \eta_k\ne\omega_1,\omega_6\}.$$
\end{enumerate}
In the following, for each $f\in\{g,h\}$, we consider only the case for $(\eta_1,\dots,\eta_k)\in B_f^{(k)}$.

\section{Basic properties} \label{sec: basic}
Now, we explain some basic properties of our target. 

\subsection{Relationship with known special values}
First, we describe the relationship with known special values.

\exm \label{ex: example}
Take $f\in\{g,h\}$ and $k,k_1,\dots,k_d\in\bZge{1}$ with $k_d\ge2$.
\begin{enumerate}
  \item By iterated integral representation of multiple zeta values (\cref{eq: MZV}), it holds $$\mathrm{I}_f(\omega_1\omega_0^{k_1-1}\cdots\omega_1\omega_0^{k_d-1})=\zeta(k_1,\dots,k_d).$$
  \item By substituting the integral variable $x=\frac{2\lambda}{1+\lambda^2}$, we have
  \begin{align}
    \mathrm{I}_g(\omega_2\omega_3^{k_1-1}\cdots\omega_2\omega_3^{k_d-1})=\int_0^1\frac{2\diff\lambda}{1+\lambda^2}\left(\frac{\diff \lambda}{\lambda}\right)^{k_1-1}\cdots\frac{2\diff\lambda}{1+\lambda^2}\left(\frac{\diff \lambda}{\lambda}\right)^{k_d-1}.
  \end{align}
  This matches the iterated integral representation of multiple $\widetilde{T}$-values
  \begin{align}
    \widetilde{T}(k_1,\dots,k_n):=2^n\sum_{\substack{0<m_1<\cdots<m_n\\m_j\equiv j \text{ mod }2}}\frac{(-1)^{(m_n-n)/2}}{m_1^{k_1}\cdots m_n^{k_n}},
  \end{align}
  which is defined by Kaneko--Tsumura \cite[p.4, Proposition 2.1]{KaTs22}. Then, it holds
  \begin{align}
    \mathrm{I}_g(\omega_2\omega_3^{k_1-1}\cdots\omega_2\omega_3^{k_d-1})=\widetilde{T}(k_1,\dots,k_n).
  \end{align}
  \item In our previous work (\cite[p.9, Proposition 3.2.1]{Ot23}), we showed
  \begin{align}
    \mathrm{I}_g(\omega_2\omega_0^{k-1})\in\Span_\bQ\left\{\pi^l\cdot(\log2)^{l_1}\cdot\prod_{j\ge3\text{: odd}}\zeta(j)^{l_j}~\middle|~\begin{array}{c}l\in\bZge{1}\text{: odd},~l_j\in\bZge{0}~(j\ge1\text{: odd})\\l+\sum_{j\ge1\text{: odd}}l_j=k\end{array}\right\}.
  \end{align}
\end{enumerate}
\eexm

\subsection{Explicit formula}

Here, we summarize formulas given by the definition and the basic properties of iterated integrals.

\propo[{\cite[p.220, p233]{Ch73}}] \label{prop: basic}
  For a smooth manifold $M$, $1$-forms $\eta_1,\dots,\eta_k,\eta_{k+1},\dots,\eta_{k+l}$ on $M$, and piecewise smooth path $\gamma,\gamma_1,\gamma_2:[0,1]\rightarrow M$ with $\gamma_1(1)=\gamma_2(0)$, the following properties hold.
  \begin{enumerate}
    \item (path connection) $$\int_{\gamma_1\gamma_2}\eta_1\cdots\eta_k=\sum_{j=0}^k\int_{\gamma_1}\eta_1\cdots\eta_j\int_{\gamma_2}\eta_{j+1}\cdots\eta_k.$$
    \item (path reversal)$$\int_{\gamma^{-1}}\eta_1\cdots\eta_k=(-1)^k\int_\gamma\eta_k\cdots\eta_1.$$
    \item (shuffle relation)$$\int_\gamma\eta_1\cdots\eta_k\int_\gamma\eta_{k+1}\cdots\eta_{k+l}=\sum_{\sigma\in S_{k,l}}\int_\gamma\eta_{\sigma(1)}\cdots\eta_{\sigma(k+l)}.$$
    \item (homotopy invariance) If each of $\eta_1,\dots,\eta_k$ is closed and satisfy $\eta_j\wedge\eta_{j+1}=0~(j=1,\dots,k-1)$, the value $\int_\gamma\eta_1\cdots\eta_k$ only depends on the homotopy class of $\gamma$. 
  \end{enumerate}
\epropo

By \cref{prop: basic}-(3), we have the following proposition.
\propo \label{prop: shuffle}
  For each $f\in\{g,h\}$ and $(\varphi_1,\dots,\varphi_k)\in B_f^{(k)}, (\varphi_{k+1},\dots,\varphi_{k+l})\in B_f^{(l)}$, it holds
  \begin{align}
    \mathrm{I}_f(\varphi_1\cdots\varphi_k)\mathrm{I}_f(\varphi_{k+1}\cdots\varphi_{k+l})=\sum_{\sigma\in S_{k,l}}\mathrm{I}_f(\varphi_{\sigma(1)}\cdots\varphi_{\sigma(k+l)}).
  \end{align}
  In particular, it holds that
  \begin{enumerate}
    \item for each $f\in\{g,h\}$, $k\in\bZge{0}$, and $(\varphi)\in B_f^{(1)}$,
    \begin{align}
      \mathrm{I}_f(\varphi^k)=\frac{1}{k!}\mathrm{I}_f(\varphi)^k;
    \end{align}
    \item for each $f\in\{g,h\}$, $k\in\bZge{0}$, $j=0,\dots,k-2$, and $\varphi,\eta\in B_f$ with $(\{\varphi\}^j,\eta,\{\varphi\}^{k-j-1})\in B_f^{(k)}$,
    \begin{align}
      \mathrm{I}_f(\varphi^j\eta\varphi^{k-j-1})=\frac{1}{k-j-1}\left\{\mathrm{I}_f(\varphi^j\eta\varphi^{k-j-2})\mathrm{I}_f(\varphi)-(j+1)\mathrm{I}_f(\varphi^{j+1}\eta\varphi^{k-j-2})\right\}.
    \end{align}
  \end{enumerate}

\epropo

Now, we give the iterated sum representation of our target.

\propo \label{prop: sum-rep}
  For each $f\in\{g,h\}$ and $(\varphi_1,\dots,\varphi_k)\in B_f^{(k)}$, it holds
  \begin{align} \label{eq: sum-rep}
    \mathrm{I}_f(\varphi_1\cdots\varphi_k)=\sum_{0< n_1\le\dots\le n_k}\prod_{j=1}^k\frac{c(n_j-n_{j-1},\varphi_j)}{n_j},
  \end{align}
  where $n_0=0$ and 
  \begin{align}
    c(n,\omega_0)&=\delta_{n,0},  & c(n,\omega_1)&=1-\delta_{n,0}\\
    c(n,\omega_2)&=\begin{cases}
      \frac{1}{2^{n-1}}\binom{n-1}{\frac{n-1}{2}}, & n \text{ is odd},\\
      0, & n \text{ is even},
    \end{cases} &
    c(n,\omega_3)&=\begin{cases}
      \frac{1}{2^{n}}\binom{n}{\frac{n}{2}}, & n \text{ is even},\\
      0, & n \text{ is odd},
    \end{cases}\\
    c(n,\omega_4)&=\begin{cases}
      \left(\frac{\sqrt{3}}{4}\right)^{n-1}\binom{n-1}{\frac{n-1}{2}}, & n \text{ is odd},\\
      0, & n \text{ is even},
    \end{cases} &
    c(n,\omega_5)&=\begin{cases}
      \left(\frac{\sqrt{3}}{4}\right)^{n}\binom{n}{\frac{n}{2}}, & n \text{ is even},\\
      0, & n \text{ is odd},
    \end{cases}\\
    c(n,\omega_6)&=\frac{1}{2}\sum_{m=0}^{\lfloor\frac{n-1}{2}\rfloor}\left(\frac{3}{16}\right)^m\binom{2m}{m}& &
  \end{align}
  for $n\in\bZge{0}$.
\epropo
\pf
  It is known that
  \begin{align}
    \int_\dch\varphi_1\cdots\varphi_k=\sum_{0< n_1\le\dots\le n_k}\prod_{j=1}^k\frac{c(n_j-n_{j-1},\varphi_j)}{n_j},
  \end{align}
  where $\varphi_j=\sum_{n=0}^\infty c(n,\varphi_j)x^{n-1}\diff x$ (see \cite[p.8, Proposition 3.1.3]{Ot23}). It implies \cref{eq: sum-rep} since $\omega_s=\sum_{n=0}^\infty c(n,\omega_s)x^{n-1}\diff x$ for each $s=0,1,\dots,6$.
\epf

Now, we investigate a special case of which our target can be written in terms of  well-known special values as logarithms, Riemann zeta values, and Dirichlet $L$-values.

\propo \label{prop: special-case1}
  For each positive integer $m\ge2$, it holds
  \begin{align}
    \mathrm{I}_g(\omega_2^{2m-1}\omega_0)&=-2^{1-2m}\sum_{l=1}^{m-1}\frac{(-1)^l}{(2m-2l-1)!}(2^{-2l}-1)\zeta(2l+1)\pi^{2m-2l-1}\\
    &~~~~+\frac{1}{(2m-1)!}\log2\left(\frac{\pi}{2}\right)^{2m-1}.\label{eq: prop33-1}
  \end{align}
  Also, for each positive integer $m\ge1$, it holds
  \begin{align}
    \mathrm{I}_g(\omega_2^{2m}\omega_0)&=-2^{-2m}\sum_{l=1}^{m-1}\frac{(-1)^l}{(2m-2l)!}(2^{-2l}-1)\zeta(2l+1)\pi^{2m-2l}\\
    &~~~~+(-1)^m2^{-2m}(2-2^{-2m})\zeta(2m+1)+\frac{1}{(2m)!}\log2\left(\frac{\pi}{2}\right)^{2m}. \label{eq: prop33-2}
  \end{align}
  In particular,
  \begin{align}
    \mathrm{I}_g(\omega_2^{k-1}\omega_0)\in\mathrm{Span}_\bQ\left\{\log2\cdot\pi^{k-1}, \zeta(2l+1)\pi^{k-2l-1}~\middle|~l=1,\dots,\left\lfloor\frac{k-1}{2}\right\rfloor\right\}.
  \end{align}
  Also, for $k\ge2$ and $1\le j\le k-2$, it holds
  \begin{align}
    \mathrm{I}_g(\omega_2^j\omega_0\omega_2^{k-j-1})\in\mathrm{Span}_\bQ\left\{\zeta(2l+1)\pi^{k-2l-1}~\middle|~l=1,\dots,\left\lfloor\frac{k-1}{2}\right\rfloor\right\}.
  \end{align}  
\epropo

To prove this explicit formula, we need the following lemma. 

\lem\label{lem: special-case1}
  If $z_0\in\bC$ satisfies $|z_0|=1$ and $-\frac{\pi}{2}\le\arg(z_0)\le\frac{\pi}{2}$, then 
  \begin{align}
    \int_{\gamma}\frac{\log(1-z)}{z}\diff z&=-\Li_2(z_0)+\frac{\pi\sqrt{-1}}{2}\log z_0+\frac{\pi^2}{6}; \label{eq: lem1-1}\\
    \int_{\gamma}\frac{\log(1+z)}{z}\diff z&=-\Li_2(-z_0)-\frac{\pi^2}{12}, \label{eq: lem1-2}
  \end{align}
  where $\gamma=C_0(1,z_0)$ and 
  \begin{align}
    \Li_k(z_0):=-\mathrm{I}_{\dch_{0,z_0}}(0;1,\{0\}^{k-1};z_0)=\sum_{n=1}^\infty\frac{z_0^n}{n^k}.
  \end{align}
\elem
\pf
  First we show the \cref{eq: lem1-1}. We have
  \begin{align}
    \int_\gamma\frac{\log(1-z)}{z}\diff z&=-\int_\gamma\diff z_2 \int_{\dch_{0,z_2}}\diff z_1 \frac{1}{1-z_1}\frac{1}{z_2}\\
    &=-\int_\gamma\diff z_2\int_{\gamma_1\gamma_2\gamma_3}\diff z_1\frac{1}{1-z_1}\frac{1}{z_2},
  \end{align}
  where $\gamma_1:=\dch_{0,1-\varepsilon}$, $\gamma_2:=C_1(1-\varepsilon,1+\varepsilon\sqrt{-1})$, and $\gamma_3:=C_0(1+\varepsilon\sqrt{-1},z_2)$ for sufficiently small $\varepsilon>0$ (see FIGURE \ref{fig: integral paths1}).
  
  \begin{figure}[htbp]
    \centering
    \caption{integral paths}
    \begin{tikzpicture}[node/.style={draw, circle, font=\large, inner sep=6pt}]
      \node [coordinate] (0) at (0,0) {};
      \node [below=0.2cm of 0] {$0$};
      \node [coordinate] (1) at (5,0) {};
      \node [right=0.2cm of 1] {$1$};
      \node [coordinate] (1-e) at (4.5,0) {};
      \node [below=0.2cm of 1-e] {$1-\varepsilon$};
      \node [coordinate] (1+ei) at (5,0.5) {};
      \node [right=0.2cm of 1+ei] {};
      \node [coordinate] (z2) at (4,3) {};
      \node [right=0.2cm of z2] {$z_2$};

      \fill (1) circle [radius=2pt];      
      \fill (0) circle [radius=2pt];
      \fill (z2) circle [radius=2pt];

      \draw[->] (0) -- (2,1.5);
      \draw (2,1.5) -- (z2);
      \draw[->] (0) -- (2.25,0);
      \draw (2.25,0) -- (1-e);
      \draw[->] (1-e) arc [start angle=180, end angle=135, radius=0.5];
      \draw (1-e) arc [start angle=180, end angle=90, radius=0.5];
      \draw[->] (1+ei) arc [start angle=5.71, end angle=21.29, radius=5];
      \draw (1+ei) arc [start angle=5.71, end angle=36.87, radius=5];

      \node at (5.2,1.5) {$\gamma_3$};
      \node at (4.3,0.7) {$\gamma_2$};
      \node at (2.25,-0.3) {$\gamma_1$};
      \node at (1.5,1.7) {$\dch_{0,z_2}$};
    \end{tikzpicture}
    \label{fig: integral paths1}
  \end{figure}

  By $\varepsilon\searrow0$ and calculating the regularization of the iterated integral (\cref{rem: reg of II}), then
  \begin{align}
    &\int_{\gamma_1}\frac{\diff z_1}{1-z_1}=0,\\
    &\int_{\gamma_2}\frac{\diff z_1}{1-z_1}=-\frac{\pi\sqrt{-1}}{2}.
  \end{align}
  Also, 
  \begin{align}
    \int_\gamma\diff z_2\int_{\gamma_3}\diff z_1\frac{1}{1-z_1}\frac{1}{z_2}=\int_\gamma\frac{\diff z}{1-z}\frac{\diff z}{z}=\Li_2(z_0)-\Li_2(1)=\Li_2(z_0)-\frac{\pi^2}{6}.
  \end{align}
  Then, we have the \cref{eq: lem1-1}. Next, we show the \cref{eq: lem1-2}. We have
  \begin{align}
    \int_\gamma\frac{\log(1+z)}{z}\diff z&=\int_\gamma\diff z_2\int_{\dch_{0,z_2}}\diff z_1\frac{1}{1+z_1}\frac{1}{z_2}\\
    &=\int_\gamma\diff z_2\left(\int_{\dch\gamma_3}\diff z_1\right)\frac{1}{1+z_1}\frac{1}{z_2},
  \end{align}
  where $\gamma_3:=C_0(1,z_2)$. Here,
  \begin{align}
    \int_\dch\frac{\diff z_1}{1+z_1}=\log2.
  \end{align}
  On the other hands, it holds
  \begin{align}
    \int_\gamma\diff z_2\int_{\gamma_3}\diff z_1\frac{1}{1+z_1}\frac{1}{z_2}&=\int_\gamma\frac{\diff z}{1+z}\frac{\diff z}{z}\\
    &=-\int_{-\gamma}\frac{\diff z}{1-z}\frac{\diff z}{z}\\
    &=-\int_{\gamma_4\gamma_5\gamma_6}\frac{\diff z}{1-z}\frac{\diff z}{z},
  \end{align}
  where $-\gamma$ is defined by $(-\gamma)(t):=-\gamma(t)$, $\gamma_4:=\dch_{-1,-\varepsilon}$, $\gamma_5:=C_0\left(-\varepsilon,-\varepsilon\frac{z_0}{|z_0|}\right)$, and $\gamma_6:=\dch_{-\varepsilon\frac{z_0}{|z_0|},-z_0}$ for sufficiently small $\varepsilon>0$ (see FIGURE \ref{fig: integral paths2}). 
  
  \begin{figure}[htbp]
    \centering
    \caption{integral paths}
    \begin{tikzpicture}[node/.style={draw, circle, font=\large, inner sep=6pt}]
      \node [coordinate] (0) at (0,0) {};
      \node [above=0.2cm of 0] {$0$};
      \node [coordinate] (-1) at (-5,0) {};
      \node [above=0.2cm of -1] {$-1$};
      \node [coordinate] (-e) at (-0.5,0) {};
      \node [above=0.2cm of -e] {$-\varepsilon$};
      \node [coordinate] (et) at (-0.4,-0.3) {};
      \node [below=0.2cm of et] {};
      \node [coordinate] (-z0) at (-4,-3) {};
      \node [below=0.2cm of -z0] {$-z_0$};

      \fill (0) circle [radius=2pt];
      \fill (-1) circle [radius=2pt];
      \fill (-z0) circle [radius=2pt];

      \draw[->] (-1) -- (-2.5,0);
      \draw (-2.5,0) -- (-e);
      \draw[->] (et) -- (-2.2,-1.65);
      \draw (-2.2,-1.65) -- (-z0);
      \draw[->] (-e) arc [start angle=180, end angle=198.435, radius=0.5];
      \draw (-e) arc [start angle=180, end angle=216.87, radius=0.5];
      \draw[->] (-1) arc [start angle=180, end angle=198.435, radius=5];
      \draw (-1) arc [start angle=180, end angle=216.87, radius=5];

      \node at (-5.5,-1.5) {$-\gamma$};
      \node at (-2.5,0.3) {$\gamma_4$};
      \node at (-1.0,-0.3) {$\gamma_5$};
      \node at (-2.0,-2.0) {$\gamma_6$};
    \end{tikzpicture}
    
    \label{fig: integral paths2}
  \end{figure}
  Then, by $\varepsilon\searrow0$, \cref{prop: basic} (1), and
  \begin{align}
    \int_{\gamma_6}\frac{\diff z}{z}=\int_{\gamma_5}\frac{\diff z}{1-z}=0,
  \end{align}
  we have
  \begin{align}
    \int_{\gamma_4\gamma_5\gamma_6}\frac{\diff z}{1-z}\frac{\diff z}{z}=\int_{\gamma_4}\frac{\diff z}{1-z}\frac{\diff z}{z}+\int_{\gamma_4}\frac{\diff z}{1-z}\int_{\gamma_5}\frac{\diff z}{z}+\int_{\gamma_6}\frac{\diff z}{1-z}\frac{\diff z}{z}.
  \end{align}
  Here, each term can be calculated as
  \begin{align}
    \int_{\gamma_4}\frac{\diff z}{1-z}\frac{\diff z}{z}&=\int_{\gamma_4^{-1}}\frac{\diff z}{z}\frac{\diff z}{1-z}\\
    &=\int_{\dch_{0,-1}}\frac{\diff z}{z}\int_{\dch_{0,-1}}\frac{\diff z}{1-z}-\int_{\dch_{0,-1}}\frac{\diff z}{1-z}\frac{\diff z}{z}\\
    &=-\Li_2(-1),\\
    \int_{\gamma_4}\frac{\diff z}{1-z}\int_{\gamma_5}&\frac{\diff z}{z}=\log2\log z_0,\\
    \int_{\gamma_6}\frac{\diff z}{1-z}\frac{\diff z}{z}&=\int_{\dch_{0,-z_0}}\frac{\diff z}{1-z}\frac{\diff z}{z}=\Li_2(-z_0),
  \end{align}
  then we have the \cref{eq: lem1-2}.
\epf

\pf[of \cref{prop: special-case1}]
  By the definition and the substitution $x=\sin t$, for $k\ge2$,
  \begin{align}
    \mathrm{I}_g(\omega_2^{k-1}\omega_0)&=\int_0^1\underbrace{\frac{\diff x}{\sqrt{1-x^2}}\cdots\frac{\diff x}{\sqrt{1-x^2}}}_{k-1}\frac{\diff x}{x}\\
    &=\int_0^{\frac{\pi}{2}}\underbrace{\diff t\cdots\diff t}_{k-1}\cot t\diff t\\
    &=\frac{1}{(k-1)!}\int_0^{\frac{\pi}{2}}t^{k-1}\cot t\diff t\\
    &=\frac{1}{(k-1)!}\left\{\left[t^{k-1}\log\sin t\right]_0^{\frac{\pi}{2}}-(k-1)\int_0^{\frac{\pi}{2}}t^{k-2}\log\sin t\diff t\right\}\\
    &=-\frac{1}{(k-2)!}\int_0^{\frac{\pi}{2}}t^{k-2}\log\sin t\diff t. \label{eq: prop33-3}
  \end{align}
  Then, if we put $z=e^{\sqrt{-1}t}$,
  \begin{align}
    \int_0^{\frac{\pi}{2}}t^{k-2}\log\sin t\diff t=\frac{1}{\sqrt{-1}^{k-1}}\int_\gamma\log^{k-2}z\log\frac{1-z^2}{-2\sqrt{-1}z}\frac{\diff z}{z},
  \end{align}
  where $\gamma=C_0(1,\sqrt{-1})$. Since $0\le\arg(z)\le\frac{\pi}{2}$ for any $z\in\gamma([0,1])$, we have
  \begin{align}
    -\frac{\pi}{2}\le\arg(1-z)\le-\frac{\pi}{4}, && 0\le\arg(1+z)\le\frac{\pi}{4}, && -\frac{\pi}{2}\le\arg(-2\sqrt{-1}z)\le0.
  \end{align}
  Then, by \cref{lem: special-case1} and the partial integration, we have
  \begin{align}
    &\int_\gamma\log^{k-2}z\log\frac{1-z^2}{-2\sqrt{-1}z}\frac{\diff z}{z}\\
    &=\left[\log^{k-2}z(-\Li_2(z)+\frac{\pi\sqrt{-1}}{2}\log z-\Li_2(-z)-\log2\log z-\frac{1}{2}\log^2z)\right]_1^{\sqrt{-1}}\\
    &~~~~-(k-2)\int_1^{\sqrt{-1}}\log^{k-3}z(-\Li_2(z)+\frac{\pi\sqrt{-1}}{2}\log z-\Li_2(-z)-\log2\log z-\frac{1}{2}\log^2z)\frac{\diff z}{z}.
  \end{align}
  Furthermore, by repeating the partial integration, we have
  \begin{align}
    &\int_\gamma\log^{k-2}z\log\frac{1-z^2}{-2\sqrt{-1}z}\frac{\diff z}{z}\\
    &=\sum_{j=2}^k(-1)^{j-1}\frac{(k-2)!}{(k-j)!}\left\{\left(\frac{\pi\sqrt{-1}}{2}\right)^{k-j}\left(\Li_j(\sqrt{-1})+\Li_j(-\sqrt{-1})\right)\right.\\
    &~~~~~~~~~~~~~~~~~~~~~~~~~~~~~~~~~~~~\left.-\frac{j-1}{j!}\left(\frac{\pi\sqrt{-1}}{2}\right)^k+\frac{1}{(j-1)!}\log2\left(\frac{\pi\sqrt{-1}}{2}\right)^{k-1}\right\}\\
    &~~~~+(-1)^k(k-2)!\left(\Li_k(1)+\Li_k(-1)\right).
  \end{align}
  Since $\mathrm{I}_g(\omega_2^{k-1}\omega_0)$ is a real number, we consider only the real part of the \cref{eq: prop33-3}. Also, it is known that
  \begin{align}
    &\Li_j(\sqrt{-1})+\Li_j(-\sqrt{-1})=2^{1-j}(2^{1-j}-1)\zeta(j);\\
    &\Li_j(1)+\Li_j(-1)=2^{1-j}\zeta(j)
  \end{align} 
  then we have
  \begin{align}
    \mathrm{I}_g(\omega_2^{k-1}\omega_0)&=-\frac{1}{(k-2)!}\left\{\sum_{\substack{j=2\\j: \text{odd}}}^k(-1)^{\frac{1-j}{2}}\frac{(k-2)!}{(k-j)!}\left(\frac{\pi}{2}\right)^{k-j}2^{1-j}(2^{1-j}-1)\zeta(j)\right.\\
    &~~~~~~~~~~~~~~~~~~~~~\left.+\sum_{j=2}^k\frac{(k-2)!}{(k-j)!(j-1)!}(-1)^{j-1}\log2\left(\frac{\pi}{2}\right)^{k-1}\right.\\
    &~~~~~~~~~~~~~~~~~~~~~\left.+\frac{1-(-1)^k}{2}(-1)^{\frac{1-k}{2}}(k-2)!2^{1-k}\zeta(k)\right\}\\
    &=-\sum_{\substack{j=2\\j: \text{odd}}}^k(-1)^{\frac{1-j}{2}}\frac{1}{(k-j)!}\left(\frac{\pi}{2}\right)^{k-j}2^{1-j}(2^{1-j}-1)\zeta(j)\\
    &~~~~+\frac{1}{(k-1)!}\log2\left(\frac{\pi}{2}\right)^{k-1}+\frac{1-(-1)^k}{2}(-1)^{\frac{1-k}{2}}2^{1-k}\zeta(k). \label{eq: prop33-4}
  \end{align}
  The \cref{eq: prop33-1,eq: prop33-2} are shown by setting $k=2m$ and $k=2m+1$. Furthermore, by \cref{prop: shuffle},
  \begin{align}
    \mathrm{I}_g(\omega_2^{k-2}\omega_0\omega_2)=\mathrm{I}_g(\omega_2^{k-2}\omega_0)\mathrm{I}_g(\omega_2)-(k-1)\mathrm{I}_g(\omega_2^{k-1}\omega_0).
  \end{align}
  Thus, by the \cref{eq: prop33-4} and $\mathrm{I}_g(\omega_2)=\frac{\pi\sqrt{-1}}{2}$, the terms containing $\log2$ vanish and we have
  \begin{align}
    \mathrm{I}_g(\omega_2^{k-2}\omega_0\omega_2)\in\mathrm{Span}_\bQ\left\{\zeta(2l+1)\pi^{k-2l-1}~\middle|~l=1,\dots,\left\lfloor\frac{k-1}{2}\right\rfloor\right\}.
  \end{align}
  Also, by repeating the equation
  \begin{align}
    \mathrm{I}_f(\varphi^j\eta\varphi^{k-j-1})=\frac{1}{k-j-1}\left\{\mathrm{I}_f(\varphi^j\eta\varphi^{k-j-2})\mathrm{I}_f(\varphi)-(j+1)\mathrm{I}_f(\varphi^{j+1}\eta\varphi^{k-j-2})\right\},
  \end{align}
  we have
  \begin{align}
    \mathrm{I}_g(\omega_2^j\omega_0\omega_2^{k-j-1})\in\mathrm{Span}_\bQ\left\{\zeta(2l+1)\pi^{k-2l-1}~\middle|~l=1,\dots,\left\lfloor\frac{k-1}{2}\right\rfloor\right\}
  \end{align}
  if $\le j\le k-2$.
\epf

\propo \label{prop: special-case2}
  For each positive integer $m\ge2$, it holds
  \begin{align}
    \mathrm{I}_h(\omega_4^{2m-1}\omega_0)&=\frac{1}{\sqrt{3}^{2m-1}}\frac{1}{(2m-1)!}\left(\frac{\pi}{\sqrt{3}}\right)^{2m-1}\frac{1}{2}\log3\\
    &~~~~+\frac{1}{3^{2m-1}}\sum_{l=1}^m\frac{(-1)^{1-l}}{(2m-2l)!}\left(\frac{\pi}{\sqrt{3}}\right)^{2m-2l}3^l2^{-2l}L(2l,\chi_{-3})\\
    &~~~~-\frac{1}{3^{2m-1}}\sum_{l=1}^{m-1}\frac{(-1)^{-l}}{(2m-2l-1)!}\left(\frac{\pi}{\sqrt{3}}\right)^{2m-2l-1}3^l2^{-2l-1}(3^{-2l}-1)\zeta(2l+1). \label{eq: prop34-1}
  \end{align}
  Here, $\chi_{-3}:\bZ\rightarrow\bC$ is the Dirichlet character defined by 
  $$\chi_{-3}(n)=\begin{cases}
    1, & n\equiv1\text{ mod $3$},\\
    -1 & n\equiv-1\text{ mod $3$},\\
    0, & n\equiv0\text{ mod $3$},
  \end{cases}$$ 
  and $$L(j,\chi_{-3}):=\sum_{n\ge1}\frac{\chi_{-3}(n)}{n^j}$$ is the Dirichlet $L$-value of $\chi_{-3}$. 
  Also, for each positive integer $m\ge1$, it holds
  \begin{align}
    \mathrm{I}_h(\omega_4^{2m}\omega_0)&=\frac{1}{3^{m}}\frac{1}{(2m)!}\left(\frac{\pi}{3}\right)^{2m}\frac{1}{2}\log3+(-1)^{-m}\frac{1}{12^m}\zeta(2m+1)\\
    &~~~~+\frac{1}{3^{2m}}\sum_{l=1}^m\frac{(-1)^{1-l}}{(2m-2l+1)!}\left(\frac{\pi}{\sqrt{3}}\right)^{2m-2l+1}3^l2^{-2l}L(2l,\chi_{-3})\\
    &~~~~-\frac{1}{3^{2m}}\sum_{l=1}^m\frac{(-1)^{-l}}{(2m-2l)!}\left(\frac{\pi}{\sqrt{3}}\right)^{2m-2l}3^l2^{-2l-1}(3^{-2l}-1)\zeta(2l+1). \label{eq: prop34-2}
  \end{align}
  In particular,
  \begin{align}
    &\mathrm{I}_h(\omega_4^{k-1}\omega_0)\\
    &\in\mathrm{Span}_\bQ\left\{\log3\cdot\left(\frac{\pi}{\sqrt{3}}\right)^{k-1},\zeta(2l+1)\left(\frac{\pi}{\sqrt{3}}\right)^{k-2l-1}, L(2l',\chi_{-3})\left(\frac{\pi}{\sqrt{3}}\right)^{k-2l'}~\middle|~\begin{array}{l}l=1,\dots,\left\lfloor\frac{k-1}{2}\right\rfloor,\\ l'=1,\dots,\left\lfloor\frac{k}{2}\right\rfloor\end{array}\right\}.
  \end{align}
  Also, for $k\ge2$ and $1\le j\le k-2$, it holds
  \begin{align}
    \mathrm{I}_h(\omega_4^j\omega_0\omega_4^{k-j-1})\in\mathrm{Span}_\bQ\left\{\zeta(2l+1)\left(\frac{\pi}{\sqrt{3}}\right)^{k-2l-1}, L(2l',\chi_{-3})\left(\frac{\pi}{\sqrt{3}}\right)^{k-2l'}~\middle|~\begin{array}{l} l=1,\dots,\left\lfloor\frac{k-1}{2}\right\rfloor,\\ l'=1,\dots,\left\lfloor\frac{k}{2}\right\rfloor\end{array}\right\}.
  \end{align}  
\epropo

\pf
  By the similar calculation as in the proof of \cref{prop: special-case1}, for $k\ge2$, we have
  \begin{align}
    \mathrm{I}_h(\omega_4^{k-1}\omega_0)&=\frac{1}{\sqrt{3}^{k-1}}\frac{1}{(k-1)!}\left\{\left(\frac{\pi}{3}\right)^{k-1}\log\frac{\sqrt{3}}{2}-\int_0^{\frac{\pi}{3}}t^{k-2}\log\sin t\diff t\right\}\\
    &=\frac{1}{\sqrt{3}^{k-1}}\frac{1}{(k-1)!}\left(\frac{\pi}{3}\right)^{k-1}\log\frac{\sqrt{3}}{2}\\
    &+\frac{1}{\sqrt{-3}^{k-1}}\frac{1}{(k-1)!}\bigg\{(-1)^k(k-2)!(\Li_k(1)+\Li_k(-1))\\
    &~~~~~~~~-\sum_{j=2}^k(-1)^{j-1}\frac{(k-1)!}{(k-j)!}\left(\frac{\pi\sqrt{-1}}{3}\right)^{k-j}\bigg(\Li_j(\xi_6)+\Li_j(-\xi_6)\\
    &~~~~~~~~~~~~~~\left.\left.+\frac{1}{(j-1)!}\log2\left(\frac{\pi\sqrt{-1}}{3}\right)^{j-1}+\frac{2-3j}{2j!}\left(\frac{\pi\sqrt{-1}}{3}\right)^j\right)\right\}. \label{eq: prop34-3}
  \end{align}
  Here, for each $j\ge2$,
  \begin{align}
    \Li_j(\xi_6)+\Li_j(-\xi_6)&=\sum_{n\ge1}\frac{\xi_6^n+(-\xi_6)^n}{n^j}\\
    &=2^{1-j}\sum_{n\ge1}\frac{\xi_3^n}{n^j}\\
    &=2^{-j}(3^{1-j}-1)\zeta(j)+2^{-j}\sqrt{-3}L(j,\chi_{-3}).
  \end{align}
  Then, the real part of the \cref{eq: prop34-3} is
  \begin{align}
    &\frac{1}{\sqrt{3}^{k-1}}\frac{1}{(k-1)!}\left(\frac{\pi}{3}\right)^{k-1}\frac{1}{2}\log3+\frac{1-(-1)^k}{2}(-1)^{\frac{1-k}{2}}\frac{1}{\sqrt{3}^{k-1}}2^{1-k}\zeta(k)\\
    &~~~~+\frac{1}{\sqrt{3}^{k-2}}\sum_{j=2\text{: even}}^k\frac{(-1)^{\frac{2-j}{2}}}{(k-j)!}\left(\frac{\pi}{3}\right)^{k-j}2^{-j}L(j,\chi_{-3})\\
    &~~~~-\frac{1}{\sqrt{3}^{k-1}}\sum_{j=3\text{: odd}}^k\frac{(-1)^{\frac{1-j}{2}}}{(k-j)!}\left(\frac{\pi}{3}\right)^{k-j}2^{-j}(3^{1-j}-1)\zeta(j)
  \end{align}
  and the \cref{eq: prop34-1,eq: prop34-2} are shown by setting $k=2m$ and $k=2m+1$. The latter claims are also shown in the same way as the proof of \cref{prop: special-case1}.
\epf

\cref{ex: example}-(2), (3) and \cref{prop: special-case1} suggest that $\mathrm{I}_g$ can be related to multiple $L$-values of level 4 and the parity of the powers of $\pi$ can be related to the parity of the number of $\omega_2$. Similarly, \cref{prop: special-case2} suggests that $\mathrm{I}_h$ can be related to multiple $L$-values of level 6 and the parity of the powers of $\pi$ can be related to the parity of the number of $\omega_4$. These observations can be understood through \cref{sec: results}.

\section{Motivic Backgrounds} \label{sec: motivic}
Here, we review some properties of the motivic iterated integrals studied by Deligne, Glanois, and Goncharov.

\subsection{Motivic iterated integrals}
Let $\overline\bQ$ be an algebraic closure of $\bQ$ and $F(\subseteq\overline\bQ)$ be a number field. Then, there is a graded commutative Hopf algebra $\cA_F$ and a graded commutative $\cA_F$-comodule $\cH_F$. We investigate the motivic structure of periods in these algebras. Here, we review some properties of them and define motivic iterated integrals. 

\thm[{\cite[p.211, section 1.1]{Go05}}] \label{thm: HF}
Let $\cA_F$ be the fundamental Hopf algebra of the abelian category $\mathrm{MT}(F)$ of mixed Tate motives over $F$. Put $\cH_F:=\cA_F\otimes_\bQ\bQ[\tau]$ where $\tau$ is a formal element. Then, $\cA_F$ and $\cH_F$ satisfy following conditions.
\begin{enumerate}[(1)]
  \item There exists a $\bQ$-algebra homomorphism $\per: \cH_F\rightarrow\bC$ such that $\per(\tau)=2\pi\sqrt{-1}$.
  \item $\cA_F$ is a graded commutative Hopf algebra. Let $\cA_F^{(k)}$ be the $k$-th degree part of $\cA_F$. In particular, there exists a $\bQ$-linear homomorphism $$\Delta: \cA_F\rightarrow\cA_F\otimes_\bQ\cA_F,$$ which is called the coproduct on $\cA_F$.
  \item $\cH_F$ is a graded commutative $\cA$-comodule. Let $\cH_F^{(k)}$ be the $k$-th degree part of $\cH_F$. In particular, there exists a $\bQ$-linear homomorphism $$\Delta: \cH_F\rightarrow\cA_F\otimes_\bQ\cH_F,$$ which is called the coaction on $\cH_F$.
  \item $\cA_F^{(0)}=\cH_F^{(0)}=\bQ$.
  \item There exists a $\bQ$-algebra isomorphism $\cA_F\cong \cH_F/(\tau\cdot\cH_F)$. Let $\rho:\cH_F\rightarrow\cA_F$ be the natural projection.

\end{enumerate}
Further, these objects satisfy some reasonable properties explained later.
\ethm

\rem \label{rem: period conjecture}
The period map $\per: \cH_F\rightarrow\bC$ is conjecturally injective (called period conjecture, see \cite[p.342, Conjecture 4.118]{GiFr}). This conjecture suggests that information about the structure of the periods is not lost by the motivic lifting. 
\erem 

The motivic iterated integrals yield elements in $\cH_F$.
\thm[Motivic Iterated Integrals, {\cite[p.212, Theorem 1.1]{Go05}}] \label{thm: MII}
For each positive integer $k\in\bZge{0}$, $k+2$ points $a_0,\dots,a_{k+1}\in F$, two tangential base points $p=(a_0,u), q=(a_{k+1},v)$, and a smooth path $\gamma$ from $p$ to $q$, there exists an element $$\Im_\gamma(p;a_1,\dots,a_k;q)\in\cH_F^{(k)}$$ called the motivic iterated integral, such that $$\per(\Im_\gamma(p;a_1,\dots,a_k;q))=\mathrm{I}_\gamma(p;a_1,\dots,a_k;q).$$ Furthermore, the motivic iterated integrals satisfy following conditions.
\begin{enumerate}[(1)]
  \item The element $\Im_\gamma(p;a_1,\dots,a_k;q)$ is independent of the parametrization of $\gamma$.
  \item The element $\rho(\Im_\gamma(p;a_1,\dots,a_k;q))$ is independent of the homotopy class of $\gamma$. We put $$\Ia(p;a_1,\dots,a_k;q):=\rho(\Im_\gamma(p;a_1,\dots,a_k;q))\in\cA_F^{(k)}.$$
  \item If $k=0$, $\Im_\gamma(p;q)=1$.
  \item If $k>0$ and $p=q$, $\Im_\gamma(p;a_1,\dots,a_k;p)=0$.
  \item The equations $\Im_\dch(0;0;1)=\Im_\dch(0;1;1)=0$ hold.
  \item (path reversal): For a smooth path $\gamma$ on $\bC$, $$\Im_{\gamma^{-1}}(p;a_1,\dots,a_k;q)=(-1)^k\Im_\gamma(q;a_k,\dots,a_1;p).$$
  \item (path connection): For smooth paths $\gamma_1, \gamma_2$ on $\bC$ with $\gamma_1(1)=\gamma_2(0)$,
  $$\Im_{\gamma_1\gamma_2}(p;a_1,\dots,a_k;q)=\sum_{s=1}^k\Im_{\gamma_1}(p;a_1,\dots,a_s;\gamma_1(1))\Im_{\gamma_2}(\gamma_2(0);a_{s+1},\dots,a_k;q).$$
\end{enumerate}
\ethm

\defi \label{def: ext motivic II Q-linearly}
For an element $w=e_{a_1}\dots e_{a_k}$ of non-commutative polynomial ring $W_F:=\bQ\langle e_{a}~|~a\in F\rangle$, the motivic iterated integral $\Im_\gamma(p;w;q)\in\cH_F$ is defined by
\begin{align}
\Im_\gamma(p;w;q):=\Im_\gamma(p;a_1,\dots,a_k;q),
\end{align}
and extend to whole $W_F$ by $\bQ$-linearity.
\edefi

From here, we investigate the structures of $\cA_F$ and $\cH_F$ more precisely. For the $\bQ$-subspace $R=\mathcal{O}^\times_{F,S}\otimes_\bZ\bQ\subseteq F^\times\otimes_\bZ\bQ$ for some finite set $S$ of finite places of $F$, let $\cA_{F,R}$ be a fundamental Hopf algebra of $R$. Then, it can be embedded into $\cA_F$ (see \cite[p.279, Appendix A.4]{Go05}). We put $\cH_{F,R}:=\cA_{F,R}\otimes_\bQ\bQ[\tau]$.

\thm[{\cite[p.215, Theorem 1.3]{Go05}}] \label{thm: AFR}
Let $a_0,\dots,a_{k+1}$ be elements of $F$. If $\tI(a_{i_1},a_{i_2},a_{i_3})\in R$ for any $0\le i_1<i_2<i_3\le k+1$, then it holds
\begin{align}
  \Im_{\dch_{a_0,a_{k+1}}}(a_0;a_1,\dots,a_k;a_{k+1})\in\cH_{F,R}^{(k)}.
\end{align}
\ethm

By taking $p=0$, $q=1$ and $a_1,\dots,a_k\in\widetilde\mu_N$, we have a motivic interpretation of multiple $L$-value.
\cor
For each $N\in\bZge{1}$, $k_1,\dots,k_d\in\bZge{1}$ and $\alpha_1,\dots,\alpha_d\in\mu_N$ with $(k_d,\alpha_d)\ne(1,1)$, there exits $L^{\mathfrak{m}}\left(\begin{matrix}
  k_1,\dots,k_d\\
  \alpha_1,\dots,\alpha_d  
\end{matrix}\right)\in\cH^{(k)}_{\bQ(\xi_N),R_N}$ such that
\begin{align}
  \per\left(L^{\mathfrak{m}}\left(\begin{matrix}
    k_1,\dots,k_d\\
    \alpha_1,\dots,\alpha_d  
  \end{matrix}\right)\right)=L\left(\begin{matrix}
    k_1,\dots,k_d\\
    \alpha_1,\dots,\alpha_d
  \end{matrix}\right).
\end{align}
Here, $R_N:= \mathcal{O}^\times_{\bQ(\xi_N),S_N}\otimes_\bZ\bQ\subseteq \bQ(\xi_N)^\times\otimes_\bZ\bQ$ for $S_N=\{(1-\xi_N^i)~|~i=1,\dots,N, (i,N)=1\}$.
\ecor

The Hopf structure of $\cH_{F,R}$ is given by the following proposition.

\propo[{\cite[p.279, Appendix A.4 and A.5]{Go05}}] \label{prop: structure}
There exists a non-canonical isomorphism
\begin{align}
  \Phi: \cA_{F,R}\overset{\sim}{\longrightarrow}T\left(R\oplus\bigoplus_{n=2}^\infty K_{2n-1}(F)\right),
\end{align}
where $K_{2n-1}$ is in the degree $2n-1$ Quillen $K$-theory. In particular, there is a non-canonical isomorphism of Hopf algebras
\begin{align}
  \Phi: \cA_{F,R}\overset{\sim}{\longrightarrow}\bQ\left\langle b_j^{(s)}~\middle|~\begin{array}{c}j\in\bZge{1}\\1\le s \le s_j\end{array}\right\rangle=:\cV_{F,R},
\end{align}
where $\cV_{F,R}$ is the non-commutative polynomial ring spanned by degree $j\ge1$ elements $b_j^{s}$ and $$s_j:=\begin{cases}
  \dim_\bQ R, & j=1,\\
  r_1(F)+r_2(F), & j>1 \text{ and } j \text{: odd},\\
  r_2(F), & j \text{: even}.
\end{cases}$$
Also, there is an isomorphism
\begin{align}
  \Phi: \cH_{F,R}\overset{\sim}{\longrightarrow}\cV_{F,R}\otimes_\bQ\bQ[\tau]=:\cU_{F,R},
\end{align}
\epropo

\cor
Set $d_{k}:=\dim_\bQ\cA_{F,R}^{(k)}$. If $\dim_\bQ R<\infty$, then
\begin{align}
  \sum_{k=0}^\infty d_kt^k=\left(1-(\dim_\bQ R)\cdot t-\frac{r_2t^2+(r_1+r_2)t^3}{1-t^2}\right)^{-1},
\end{align}
where $r_1:=r_1(F)$ and $r_2:=r_2(F)$. Also, set $d_{k}':=\dim_\bQ\cH_{F,R}^{(k)}$, then 
\begin{align}
  \sum_{k=0}^\infty d_kt^k=\left(1-(\dim_\bQ R)\cdot t-\frac{r_2t^2+(r_1+r_2)t^3}{1-t^2}\right)^{-1}\cdot(1-t)^{-1}.
\end{align}

\ecor

\subsection{Specific cases}
In this subsection, we review some properties of the structure of $\cH_{F,R}$ for some specific cases used in \cref{sec: results}. We consider the two cases: (I) $(F,R)=(\bQ(\xi_4),\langle2\rangle_{F_g})$, which relates to multiple $L$-values of level 4; (II) $(F,R)=(\bQ(\xi_6),\langle2,3\rangle_{F_h})$, which relates to multiple $L$-values of level 6.

\propo[{\cite[p.46, Th\'{e}or\`{e}me 8.9]{De10}}] \label{prop: span}
The following conditions hold.
\begin{enumerate}[(I)]
  \item The $\bQ$-linear space $\cH_{\bQ(\xi_4),\langle2\rangle_{F_g}}^{(k)}$ is spanned by motivic iterated integrals formed
  \begin{align}
    \Im_\dch(0;a_1,\dots,a_k;1) ~~~~~~~~ (a_1,\dots,a_k\in\widetilde\mu_4).
  \end{align}

  \item The $\bQ$-linear space $\cH_{\bQ(\xi_6),\langle2,3\rangle_{F_h}}^{(k)}$ is spanned by motivic iterated integrals formed
  \begin{align}
    \Im_\dch(0;a_1,\dots,a_k;1) ~~~~~~~~ (a_1,\dots,a_k\in\widetilde\mu_6).
  \end{align}
\end{enumerate}
\epropo

\propo[{\cite[p.49, Th\'{e}or\`{e}me 6.8]{DeGo05}}] \label{prop: allrels}
The following conditions hold.
\begin{enumerate}[(I)]
  \item For $k\in\bZge{1}$, the elements $\Ia(0;\delta,\{0\}^{k-1};1)$ $(\delta\in\mu_4)$ are subject only to the following relations in $\cA_{\bQ(\xi_4),\langle2\rangle_{F_g}}$:
  \begin{align}
    \begin{cases}
      \Ia(0;\delta,\{0\}^{k-1};1)=2^{k-1}\sum_{\epsilon^2=\delta}\Ia(0;\epsilon,\{0\}^{k-1};1), & \forall\delta\in\mu_2,\\
      \Ia(0;\delta,\{0\}^{k-1};1)=(-1)^{k-1}\Ia(0;\delta^{-1},\{0\}^{k-1};1), & \forall\delta\in\mu_4.
    \end{cases}
  \end{align}

  \item For $k\in\bZge{1}$, the elements $\Ia(0;\delta,\{0\}^{k-1};1)$ $(\delta\in\mu_6)$ are subject only to the following relations in $\cA_{\bQ(\xi_6),\langle2,3\rangle_{F_h}}$:
  \begin{align}
    \begin{cases}
      \Ia(0;\delta,\{0\}^{k-1};1)=2^{k-1}\sum_{\varepsilon^2=\delta}\Ia(0;\varepsilon,\{0\}^{k-1};1), & \forall\delta\in\mu_3,\\
      \Ia(0;\delta,\{0\}^{k-1};1)=3^{k-1}\sum_{\varepsilon^3=\delta}\Ia(0;\varepsilon,\{0\}^{k-1};1), & \forall\delta\in\mu_2,\\
      \Ia(0;\delta,\{0\}^{k-1};1)=(-1)^{k-1}\Ia(0;\delta^{-1},\{0\}^{k-1};1) & \forall\delta\in\mu_6.
    \end{cases}
  \end{align}  
\end{enumerate}
\epropo

\propo[{\cite[p.354, Corollary 3.7]{Gl16}}] \label{prop: basis}
We put $\Delta':=\Delta-1\otimes\id: \cH_{F,R}\rightarrow\cA_{F,R}\otimes_\bQ\cH_{F,R}$. Then, the following equations hold:
\begin{align}
  \Ker\Delta'\cap\cH_{\bQ(\xi_4),\langle2\rangle_{F_g}}^{(k)}=\bQ\tau^k\oplus\bQ\Im_\dch(0;\xi_4,\{0\}^{k-1};1);\\
  \Ker\Delta'\cap\cH_{\bQ(\xi_6),\langle2,3\rangle_{F_h}}^{(k)}=\bQ\tau^k\oplus\bQ\Im_\dch(0;\xi_6,\{0\}^{k-1};1).
\end{align}
\epropo

\section{Main Results} \label{sec: results}
\subsection{Base extension} \label{subsec: baseext}
First, we consider the pullback of $B_f$ under $\phi_f:\bP^1(\bC)\rightarrow X_f$ to use \cref{sec: motivic} for giving a motivic interpretation of $\mathrm{I}_f(\eta_1\cdots\eta_k)$.

\vspace{5mm}
\noindent
\textbf{Pullback of $B_f$.}

\begin{enumerate}[(I)]
  \item \underline{$g(X,Y,Z)=X^2+Y^2-Z^2$}: for each element of $B_g=\left\{\omega_0:=\frac{\diff x}{x},\omega_1:=\frac{\diff x}{1-x},\omega_2:=\frac{\diff x}{y},\omega_3:=\frac{\diff x}{xy}\right\}$, the pullback along $\phi_g$ is 
  \begin{align} \label{eq: pullback of phig}
    \begin{split}
    \phi_g^*\omega_0&=\frac{(1-\lambda^2)\diff\lambda}{\lambda(1+\lambda^2)}=\frac{\diff\lambda}{\lambda}-\frac{\diff\lambda}{\lambda-\xi_4}-\frac{\diff\lambda}{\lambda-\xi_4^{-1}},\\
    \phi_g^*\omega_1&=\frac{2(1+\lambda)\diff\lambda}{(1-\lambda)(1+\lambda^2)}=-2\frac{\diff\lambda}{1-\lambda}+\frac{\diff\lambda}{\lambda-\xi_4}+\frac{\diff\lambda}{\lambda-\xi_4^{-1}},\\
    \phi_g^*\omega_2&=\frac{2\diff\lambda}{1+\lambda^2}=\xi_4^{-1}\left(\frac{\diff\lambda}{\lambda-\xi_4}-\frac{\diff\lambda}{\lambda-\xi_4^{-1}}\right),\\
    \phi_g^*\omega_3&=\frac{\diff\lambda}{\lambda}.
    \end{split}
  \end{align}
  The integral path is $\alpha: [0,1]\rightarrow\bP^1(\bC)$, $\alpha(t)=1-\sqrt{\frac{1-t}{1+t}}$ and it satisfies $\alpha([0,1])=[0,1]$.
  Then, each $\mathrm{I}_g(\eta_1\cdots\eta_k)$ can be written by a $\bQ(\xi_4)$-linear combination of iterated integrals on $\bP^1(\bC)\setminus\{0,1,\pm\sqrt{-1},\infty\}$ having the form
  \begin{align} \label{eq: f1 pullback}
    \mathrm{I}_{\dch}(0;a_1,\dots,a_k;1)
  \end{align}
  with $a_j\in F_g:=\bQ(\xi_4)$ for $j=1,\dots,k$. Also, for each of the terms given by the \cref{eq: f1 pullback}, it holds $\widetilde{\mathrm{I}}(a_{i_1},a_{i_2},a_{i_3})\in\langle2\rangle_{F_g}=:R_g\subseteq F_g^\times\otimes_\bZ\bQ$ if $0\le i_1<i_2<i_3\le k+1$ since $a_{i_1},a_{i_2},a_{i_3}\in\{0,1,\pm\sqrt{-1}\}$.

  \item \underline{$h(X,Y,Z)=X^2+XY+Y^2-Z^2$}: for each element of $B_h=\left\{\omega_0:=\frac{\diff x}{x},\omega_1:=\frac{\diff x}{1-x},\omega_4:=\frac{\diff x}{x+2y},\omega_5:=\frac{\diff x}{x(x+2y)},\omega_6:=\frac{\diff x}{(1-x)(x+2y)}\right\}$, the pullback along $\phi_h$ is 
  \begin{align}
    \phi_h^*\omega_0&=\frac{-(\lambda^2-2\lambda-2)\diff\lambda}{\lambda(2+\lambda)(1+\lambda+\lambda^2)}=\frac{\diff\lambda}{\lambda}+\frac{\diff\lambda}{\lambda+2}-\frac{\diff\lambda}{\lambda-\xi_3}-\frac{\diff\lambda}{\lambda-\xi_3^{-1}},\\
    \phi_h^*\omega_1&=\frac{(\lambda^2-2\lambda-2)\diff\lambda}{\lambda^3-1}=-\frac{\diff\lambda}{\lambda-1}+\frac{\diff\lambda}{\lambda-\xi_3}+\frac{\diff\lambda}{\lambda-\xi_3^{-1}},\\
    \phi_h^*\omega_4&=\frac{\diff\lambda}{1+\lambda+\lambda^2}=\frac{1}{\sqrt{-3}}\left(\frac{\diff\lambda}{\lambda-\xi_3}-\frac{\diff\lambda}{\lambda-\xi_3^{-1}}\right),\\
    \phi_h^*\omega_5&=\frac{\diff\lambda}{\lambda(\lambda+2)}=\frac{1}{2}\left(\frac{\diff\lambda}{\lambda}-\frac{\diff\lambda}{\lambda+2}\right),\\
    \phi_h^*\omega_6&=-\frac{\diff\lambda}{\lambda-1}.
  \end{align}
  The integral path is $\alpha: [0,1]\rightarrow\bP^1(\bC)$, $\alpha(t)=\frac{t-2+\sqrt{4-3t^2}}{2(1-t)}$ and it satisfies $[0,1]$.
  Then, each $\mathrm{I}_h(\eta_1\cdots\eta_k)$ can be written by a $\bQ(\xi_6)$-linear combination of iterated integrals on $\bP^1(\bC)\setminus\{0,1,-2,\xi_3^{\pm1},\infty\}$ having the form
  \begin{align} \label{eq: f2 pullback}
    \mathrm{I}_{\dch}(0;a_1,\dots,a_k;1)
  \end{align}
  with $a_j\in F_h:=\bQ(\xi_6)$ for $j=1,\dots,k$. Also, for each of the terms given by the \cref{eq: f2 pullback}, it holds $\widetilde{\mathrm{I}}(a_{i_1},a_{i_2},a_{i_3})\in\langle2,3\rangle_{F_h}=:R_h\subseteq F_h^\times\otimes_\bZ\bQ$ if $0\le i_1<i_2<i_3\le k+1$ since $a_{i_1},a_{i_2},a_{i_3}\in\{0,1,-2,\xi_3^{\pm1}\}$. 
\end{enumerate}

By the above discussion, for each $f\in\{g,h\}$ and $(\eta_1,\dots,\eta_k)\in B_f^{(k)}$, $\mathrm{I}_f(\eta_1\cdots\eta_k)$ can be written as an $F_f$-linear combination of iterated integrals having the form
\begin{align} \label{eq: II on P1C}
  \mathrm{I}_{\dch}(0;a_1,\dots,a_k;1)
\end{align}
with $a_j\in F_f$ for $j=1,\dots,k$. We recall that the $\bQ$-linear combination of the \cref{eq: II on P1C} has a motivic interpretation given by motivic iterated integrals in $\cH_{F_f,R_f}$. Then, we consider the base extension of $\cH_{F_f,R_f}$ to give motivic interpretations of iterated integrals on $X_f$.

\defi \label{def: base ext}
For a number field $F\supseteq\bQ$ and $R=\mathcal{O}^\times_{F,S}\otimes_\bZ\bQ\subseteq F^\times\otimes_\bZ\bQ$, we extend the following notations to $\cH_{F,R}\otimes_\bQ F$. For $u\otimes c, u_1\otimes c_1, u_2\otimes c_2\in\cH_{F,R}\otimes_\bQ F$,
\begin{itemize}
  \item coaction: $\Delta: \cH_{F,R}\otimes_\bQ F\rightarrow \cA_{F,R}\otimes_\bQ\cH_{F,R}\otimes_\bQ F,~~\Delta(u\otimes c):=\Delta(u)\otimes c$;
  \item product: $(\cH_{F,R}\otimes_\bQ F)\otimes_\bQ(\cH_{F,R}\otimes_\bQ F)\rightarrow\cH_{F,R}\otimes_\bQ F,~~(u_1\otimes c_1)\cdot(u_2\otimes c_2):=u_1u_2\otimes c_1c_2$;
  \item period map: $\per: \cH_{F,R}\otimes_\bQ F\rightarrow\bC,~~\per(u\otimes c):=\per(u)\cdot c$;
  \item isomorphism: $\Phi: \cH_{F,R}\otimes_\bQ F\rightarrow \cU_{F,R}\otimes_\bQ F,~~\Phi(u\otimes c):=\Phi(u)\otimes c$.
\end{itemize}
\edefi

\rem \label{rem: extended period conjecture}
As same as \cite[p.19, Proposition 4.4.2]{Ot23}, if the period conjecture is true, (i.e. the period map $\per: \cH_{F,R}\rightarrow\bC$ is injective), then its extension $\per:\cH_{F,R}\otimes_\bQ F\rightarrow\bC$ is also injective. This suggests that information about the structure of the periods is not preserved by the base extension.
\erem

\defi \label{def: ext motivic II F-linearly}
Let $w=c\cdot e_{a_1}\dots e_{a_k}$ be an element of non-commutative polynomial ring $F\langle F\rangle$ with $c\in F$, $a_1,\dots,a_k\in F$ satisfying $\widetilde{\mathrm{I}}(a_p,a_q,a_r)\in R$. Then, for a smooth path $\gamma:[0,1]\rightarrow \bC$ connecting $p\in\bC$ and $q\in\bC$, the motivic iterated integral $\Im_\gamma(p;w;q)\in\cH_{F,R}\otimes_\bQ F$ is defined by
\begin{align}
\Im_\gamma(p;w;q):=\Im_\gamma(p;a_1,\dots,a_k;q)\otimes c.
\end{align}
We also extend $\Im_\gamma(p;w;q)$ $F$-linearity to all $w\in F\langle F\rangle$.
\edefi

Motivic interpretations of iterated integrals on $X_{f}$ are given as elements of $\cH_{F_f,R_f}\otimes_\bQ F_f$.

\defi \label{def: motivic interpretations}
For $f\in\{g,h\}$, we define the motivic iterated integrals on $X_f$ as follows.
\begin{enumerate}[(I)]
  \item For $(\eta_1,\dots,\eta_k)\in B_g^{(k)}$, we define
  \begin{align}
    \Im_g(\eta_1\cdots\eta_k):=\Im_\dch(0;\theta_1\cdots\theta_k;1)\in\cH_{F_g,R_g}\otimes_\bQ F_g
  \end{align}
  where $\theta_1,\cdots,\theta_k\in F_g\langle e_a~|~a\in F_g\rangle$ are defined by
  \begin{align} \label{eq: thetaj}
    \begin{split}
    \theta_j:=\begin{cases}
      e_0-e_{\xi_4}-e_{\xi_4^{-1}}, & \eta_j=\omega_0,\\
      -2e_1+e_{\xi_4}+e_{\xi_4^{-1}}, & \eta_j=\omega_1,\\
      \xi_4^{-1}(e_{\xi_4}-e_{\xi_4^{-1}}), & \eta_j=\omega_2,\\
      e_0, & \eta_j=\omega_3.
    \end{cases}
  \end{split}
  \end{align}

  \item For $(\eta_1,\dots,\eta_k)\in B_h^{(k)}$, we define
  \begin{align}
    \Im_h(\eta_1\cdots\eta_k):=\Im_\dch(0;\theta_1\cdots\theta_k;1)\in\cH_{F_h,R_h}\otimes_\bQ F_h
  \end{align}
  where $\theta_1,\cdots,\theta_k\in F_h\langle e_a~|~a\in F_h\rangle$ are defined by
  \begin{align}
    \theta_j:=\begin{cases}
      e_0+e_{-2}-e_{\xi_3}-e_{\xi_3^{-1}}, & \eta_j=\omega_0,\\
      -e_1+e_{\xi_3}+e_{\xi_3^{-1}}, & \eta_j=\omega_1,\\
      \frac{1}{\sqrt{-3}}(e_{\xi_3}-e_{\xi_3^{-1}}), & \eta_j=\omega_4,\\
      \frac{1}{2}(e_0-e_{-2}), & \eta_j=\omega_5,\\
      -e_1, & \eta_j=\omega_6.
    \end{cases}
  \end{align}

\end{enumerate}
\edefi
Then, by \cref{def: ext motivic II F-linearly} and the pullback of $B_f$, we have
\begin{align}
  \per(\Im_f(\eta_1\cdots\eta_k))=\mathrm{I}_f(\eta_1\cdots\eta_k)
\end{align}
for each $f\in\{g,h\}$ and $(\eta_1,\dots,\eta_k)\in B_f^{(k)}$. This allows us to give a motivic interpretation of the iterated integrals on $X_{f}$.

\subsection{Galois action} \label{subsec: Galoisact}
In the following, the extension $F/\bQ$ is assumed to be Galois. We consider the action of the Galois group $\Gal(F/\bQ)$ on $\cH_{F,R}\otimes_\bQ F$ to investigate the Galois fixed part of $\cH_{F,R}\otimes_\bQ F$ in detail.

Let $S$ be a finite set of finite places in $F$ and $\mathcal{O}_{F,S}$ the ring of $S$-integers of $F$. Let $\mathcal{G}_F$ be the motivic Galois group associated to the canonical fiber functor of $\mathrm{MT}(F)$ defined in \cite[Section 1.2]{DeGo05}. It is a pro-algebraic affine group scheme over $\bQ$. Let $\mathcal{O}(\mathcal{G}_F)$ be the ring of the regular functions on $\mathcal{G}_F$. 
The Galois action of $\Gal(F/\bQ)$ on $F$ induces the 
action on each object of $\mathrm{MT}(F)$. Further, this action induces the action on 
$\mathcal{O}(\mathcal{G}_F)$ and on its subalgebra $\cH_F$. This action is explicitly computed in \cite[p.212, Theorem 1.1]{Go05} for the motivic iterated integrals belonging to $\cH_F$. 
Let $\mathcal{G}^{\mathrm{dR}}_F$ be the motivic Galois group associated to the de Rham fiber functor of $\mathrm{MT}(F)$  defined in \cite[Section 2.9]{DeGo05}. It is a pro-algebraic affine group scheme over $F$. Let $\mathcal{O}(\mathcal{G}^{\mathrm{dR}}_F)$ be the ring of the regular functions on $\mathcal{G}^{\mathrm{dR}}_F$.  
Then, by \cite[Proposition 2.10]{DeGo05}, 
$\mathcal{G}^{\mathrm{dR}}_F\simeq 
\mathcal{G}_F\times_{\mathrm{Spec}\hspace{0.5mm}\bQ}\mathrm{Spec}\hspace{0.5mm}F$ and it yields 
$\mathcal{O}(\mathcal{G}^{\mathrm{dR}}_F)\simeq 
\mathcal{O}(\mathcal{G}_F)\otimes_{\bQ}F$ such that 
the Galois action of $\Gal(F/\bQ)$ on 
$\mathcal{O}(\mathcal{G}^{\mathrm{dR}}_F)$ translates into 
the diagonal action on $\mathcal{O}(\mathcal{G}_F)\otimes_{\bQ}F$. 
Now, by definition, $\cH_{F,R}$ is a $\Gal(F/\bQ)$-stable $\bQ$-subalgebra of $\mathcal{O}(\mathcal{G}_F)$ because the canonical fiber functor on $MT(F)$ is stable under the action of $\Gal(F/\bQ)$. The Galois group $\Gal(F/\bQ)$ acts diagonally on the base extension 
$\cH_{F,R}\otimes_{\bQ}F$ and this action comes from 
the Galois action of $\Gal(F/\bQ)$ on 
$\mathcal{O}(\mathcal{G}^{\mathrm{dR}}_F)$ via the above isomorphism.  
By \cref{prop: span}, each element of $\cH_{F,R}$ can be written by the motivic iterated integrals if $(F,R)=(F_g,R_g)$ or $(F,R)=(F_h,R_h)$, then we can compute the Galois action explicitly for the motivic iterated integrals 
with coefficients in $F$ as follows.

\propo \label{prop: galois}
If $(F,R)=(F_f,R_f)$ for $f\in\{g,h\}$, each element $\sigma$ of  $\Gal(F/\bQ)$ acts on the motivic iterated integrals 
with coefficients in $F$ as
\begin{align}
\sigma(\Im_\dch(0;a_1,\dots,a_k;1)\otimes c)=\Im_\dch(0;{}^\sigma{a_1},\dots,{}^\sigma{a_k};1)\otimes{}^\sigma{c},
\ c\in F.
\end{align}
\epropo

\pf
It follows from the definition of the action and \cite[p.212, Theorem 1.1]{Go05}.
\epf

Let us look back at \cref{def: motivic interpretations}. Recall that we put $G_f:=\Gal(F_f/\bQ)$.

\thm[\cref{thm: 1}]
For each $f\in\{g,h\}$ and $(\eta_1,\dots,\eta_k)\in B_f^{(k)}$, there exists an element $$\Im_f(\eta_1\cdots\eta_k)\in\left(\cH_{F_f,R_f}^{(k)}\otimes_\bQ F_f\right)^{G_f}$$
such that $\per(\Im_f(\eta_1\cdots\eta_k))=\mathrm{I}_f(\eta_1\cdots\eta_k)$.
\ethm

\pf
By definition $\Im_f(\eta_1\cdots\eta_k):=\Im_\dch(0;\theta_1\cdots\theta_k;1)\in\cH_{F_f,R_f}\otimes_\bQ F_f$ and since each $\theta_p$ is $G_f$-invariant, $\Im_f(\eta_1\cdots\eta_k)$ is also $G_f$-invariant and we have
$$\Im_f(\eta_1\cdots\eta_k)\in\left(\cH_{F_f,R_f}^{(k)}\otimes_\bQ F_f\right)^{G_f}.$$
\epf

Then, it is sufficient to investigate the $G_f$-invariant part $\left(\cH_{F_f,R_f}\otimes_\bQ F_f\right)^{G_f}$ of $\cH_{F_f,R_f}\otimes_\bQ F_f$ to discuss the motivic structure of iterated integrals on $X_{f}$.

\subsection{The structure of $\left(\cH_{F_f,R_f}\otimes_\bQ F_f\right)^{G_f}$}
Now, we investigate the structure of $\left(\cH_{F_f,R_f}\otimes_\bQ F_f\right)^{G_f}$ and reduce to the original iterated integrals $\mathrm{I}_f(\eta_1\cdots\eta_k)$.

First, in preparation, we write down \cref{prop: structure} explicitly for each case of $\cH_{F_f,R_f}$.
\begin{enumerate}[(I)]
  \item if $f=g$, we have
  \begin{align}
    \cU_{F_g,R_g}=\bQ\langle b_j~|~j\in\bZge{1}\rangle\otimes_\bQ\bQ[\tau].
  \end{align}
  In particular, the degree $k$ part $\cU_{F_g,R_g}^{(k)}$ is
  \begin{align}
    \bigoplus_{r\ge0}\Span_\bQ\left\{ b_{j_1}\cdots b_{j_r}\cdot\tau^l~\middle|~\begin{array}{c}j_1,\dots,j_r\in\bZge{1},l\in\bZge{0},\\\sum_{i=1}^r j_i+l=k\end{array}\right\}.
  \end{align}

  \item if $f=h$, we have
  \begin{align}
    \cU_{F_{h},R_{h}}=\bQ\left\langle b_j^{(s)}~\middle|~j\in\bZge{1},s\in\begin{cases}
      \{1,2\}, & j=1\\
      \{1\}, & j>1
    \end{cases}\right\rangle\otimes_\bQ\bQ[\tau].
  \end{align}
  In particular, the degree $k$ part $\cU_{F_h,R_h}^{(k)}$ is
  \begin{align}
    \bigoplus_{r\ge0}\Span_\bQ\left\{ b_{j_1}^{(s_1)}\cdots b_{j_r}^{(s_r)}\cdot\tau^l~\middle|~\begin{array}{c}j_1,\dots,j_r\in\bZge{1},l\in\bZge{0},\\s_i\in\begin{cases}\{1,2\}, & j_i=1\\ \{1\}, & j_i>1\end{cases},~\sum_{i=1}^r j_i+l=k\end{array}\right\}.
  \end{align}

\end{enumerate}

Here, we fix one isomorphism given in \cref{prop: structure}. This allows us to write down the Galois action on $\cU_{F_f,R_f}$ induced by $\widetilde\Phi:\cH_{F_f,R_f}\rightarrow\cU_{F_f,R_f}$ concisely.

\propo \label{prop: normalize}
For each $f\in\{g,h\}$, we can normalize a non-canonical isomorphism $\Phi: \cH_{F_f,R_f}\rightarrow\cU_{F_f,R_f}$ in \cref{prop: structure} as follows.
\begin{enumerate}[(I)]
  \item If $f=g$, the normalized isomorphism $\widetilde\Phi: \cH_{F_g,R_g}\rightarrow\cU_{F_g,R_g}$ is defined by
  \begin{align}
    \widetilde\Phi(u_k):=b_k,
  \end{align}
  where
  \begin{align}
    u_k:=\Im_\dch(0;\xi_4,\{0\}^{k-1};1)-(-1)^k\Im_\dch(0;\xi_4^{-1},\{0\}^{k-1};1)\in\cH_{F_g,R_g}.
  \end{align}

  \item If $f=h$, the normalized isomorphism $\widetilde\Phi: \cH_{F_h,R_h}\rightarrow\cU_{F_h,R_h}$ is defined by
  \begin{align}
    &\widetilde\Phi(v_k):=b_k^{(1)},\\
    &\widetilde\Phi(v):=b_1^{(2)},
  \end{align}
  where
  \begin{align}
    &v_k:=\Im_\dch(0;\xi_3,\{0\}^{k-1};1)-(-1)^k\Im_\dch(0;\xi_3^{-1},\{0\}^{k-1};1)\in\cH_{F_h,R_h},\\
    &v:=\Im_\dch(0;-1;1)\in\cH_{F_h,R_h}.
  \end{align}
\end{enumerate}
\epropo

\pf
(I) is shown in \cite[p.21, Proposition 5.1.1]{Ot23}. Then, we only show (II). It is sufficient to prove that each of $\{v_1,v,\tau\}$ and $\{v_k,\tau^k\}~(k>1)$ is $\bQ$-linear independent. First, by \cref{prop: allrels}, we have
\begin{align}
  \Im_\dch(0;\xi_3,\{0\}^{k-1};1)+(-1)^k\Im_\dch(0;\xi_3^{-1},\{0\}^{k-1};1)\in\tau^k\bQ
\end{align}
for any $k\ge1$. On the other hands, by \cref{prop: basis}, each of $\{v,\tau\}$ and $\{\Im_\dch(0;\xi_3^{-1},\{0\}^{k-1};1), \tau^k\}$ is $\bQ$-linear independent. Thus, $\{v_k,\tau^k\}$ is also $\bQ$-linear independent. Furthermore, looking back at \cref{prop: allrels}, there are no $\bQ$-linear relations between $v$ and $v_1$. Then, $\{v,v_1,\tau\}$ is $\bQ$-linear independent.

\epf

Now, we calculate the action of the Galois group $G_f$ on $\cU_{F_f,R_f}$ induced by the isomorphisms $\widetilde\Phi: \cH_{F_f,R_f}\rightarrow\cU_{F_f,R_f}$. First, we calculate the action on each of the generators $b_j^{(s)}$ of $\cU_{F_f,R_f}$. Recall that $G_f:=\Gal(F_f/\bQ)$.

\begin{enumerate}[(I)]
  \item If $f=g$, $G_g=\{1,\sigma\}$ where $\sigma$ is defined by $\sigma(\xi_4)=\xi_4^{-1}$. By definition of $u_k$, we have
  \begin{align}
    {}^\sigma u_k=(-1)^{k+1}u_k.
  \end{align}
  Then, we have ${}^\sigma b_k=(-1)^{k+1}b_k$. Also, we have ${}^\sigma\tau=-\tau$.

  \item If $f=h$, $G_h=\{1,\sigma\}$ where $\sigma$ is defined by $\sigma(\sqrt{-3})=-\sqrt{-3}$. By definition of $v$ and $v_k$, we have
  \begin{align}
    &{}^\sigma v=v,\\
    &{}^\sigma v_k=(-1)^{k+1}v_k.
  \end{align}
  Then, we have ${}^\sigma b_k^{(1)}=b_k^{(1)}~(k\ge1)$ and ${}^\sigma b_1^{(2)}=b_1^{(2)}$. Also, we have ${}^\sigma\tau=-\tau$.

\end{enumerate}

The Galois action on $\cU_{F,R}$ is commutative with the non-commutative product of $\cU_{F,R}$. That fact is shown in the case of $(F,R)=(\bQ(\xi_4),\langle2\rangle_{F_g})$ in {\cite[p.23, Lemma 5.2.1]{Ot23}} and we can show similarly if $(F,R)=(\bQ(\xi_6),\langle2,3\rangle_{F_h})$.

\lem \label{lem: product and action}
For the generators $b_{j_p}^{(s_p)}~(p=1,\dots,r)$ of $\cU_{F_f,R_f}$, a non-negative integer $l\ge0$, and $\sigma\in\Gal(F_f/\bQ)$ the following equation holds:
\begin{align}
{}^\sigma(b_{j_1}^{(s_1)}\dots b_{j_r}^{(s_r)}\cdot\tau^l)={}^\sigma b_{j_1}^{(s_1)}\dots {}^\sigma b_{j_r}^{(s_r)}\cdot{}^\sigma(\tau^l).
\end{align}
\elem

Then, we can identify the $G_f$-invariant part $\left(\cU_{F_f,R_f}\otimes_\bQ F_f\right)^{G_f}$ of $\cU_{F_f,R_f}\otimes_\bQ F_f$.

\propo \label{prop: direct decomposition}
The following conditions hold.
\begin{enumerate}[(I)]
  \item If $f=g$,
  \begin{align}
    \left(\cU_{F_g,R_g}^{(k)}\otimes_\bQ F_g\right)^{G_g}&=\bigoplus_{\substack{r\ge0\\k-r\text{: even}}}\Span_\bQ\left\{ b_{j_1}\cdots b_{j_r}\cdot\tau^l\otimes1~\middle|~\begin{array}{c}j_1,\dots,j_r\in\bZge{1},l\in\bZge{0},\\\sum_{i=1}^r j_i+l=k\end{array}\right\}\\
    &\oplus\bigoplus_{\substack{r\ge0\\k-r\text{: odd}}}\Span_\bQ\left\{ b_{j_1}\cdots b_{j_r}\cdot\tau^l\otimes\sqrt{-1}~\middle|~\begin{array}{c}j_1,\dots,j_r\in\bZge{1},l\in\bZge{0},\\\sum_{i=1}^r j_i+l=k\end{array}\right\}.
  \end{align}
  In particular, it holds $$\dim_\bQ\left(\cU_{F_g,R_g}^{(k)}\otimes_\bQ F_g\right)^{G_g}=D^{(k)}_g,$$ where $\sum_{k=0}^\infty D_{g}^{(k)}t^k=\frac{1}{1-2t}$.

  \item If $f=h$,
  \begin{align}
    &\left(\cU_{F_h,R_h}^{(k)}\otimes_\bQ F_h\right)^{G_h}\\
    &=\bigoplus_{\substack{r\ge0\\k-r\text{: even}}}\Span_\bQ\left\{ b_{j_1}^{(s_1)}\cdots b_{j_r}^{(s_r)}\cdot\tau^l\otimes1~\middle|~\begin{array}{c}j_1,\dots,j_r\in\bZge{1},l\in\bZge{0},\\s_i\in\begin{cases}\{1,2\}, & j_i=1\\ \{1\}, & j_i>1\end{cases},~\sum_{i=1}^r j_i+l=k\end{array}\right\}\\
    &\oplus\bigoplus_{\substack{r\ge0\\k-r\text{: odd}}}\Span_\bQ\left\{ b_{j_1}^{(s_1)}\cdots b_{j_r}^{(s_r)}\cdot\tau^l\otimes\sqrt{-3}~\middle|~\begin{array}{c}j_1,\dots,j_r\in\bZge{1},l\in\bZge{0},\\s_i\in\begin{cases}\{1,2\}, & j_i=1\\ \{1\}, & j_i>1\end{cases},~\sum_{i=1}^r j_i+l=k\end{array}\right\}.
  \end{align}
  In particular, it holds $$\dim_\bQ\left(\cU_{F_h,R_h}^{(k)}\otimes_\bQ F_d\right)^{G_h}=D_h^{(k)},$$ where $\sum_{k=0}^\infty D_h^{(k)}t^k=\frac{1}{1-3t+t^2}$.
\end{enumerate}
\epropo

Here, we reduce them to $\mathrm{I}_f(\eta_1\cdots\eta_k)\in\bC$ by the period map $\per: \cH_{F_f,R_f}\otimes_\bQ F_f\rightarrow\bC$. Since $Z_f\subseteq\per\left(\left(\cH_{F_f,R_f}\otimes_\bQ F_f\right)^{G_f}\right)$, we have the following corollary.

\cor[\cref{cor: dimension}] \label{cor: dim2}
For $f\in\{g,h\}$, it holds
\begin{align}
  \dim_\bQ Z_{f}^{(k)}\le D_f^{(k)}.
\end{align}
\ecor

\begin{table}[h]
  \label{table: dimension}
  \centering
  \caption{dimension evaluation}
   \begin{tabular}{|c|c|c|c|c|c|c|}
    \hline
    $k$ & $0$ & $1$ & $2$ & $3$ & $4$ & $5$\\
    \hline \hline
    $\# B_g^{(k)}$ & & $1$ & $6$ & $24$ & $96$ & $384$\\
    \hline
    $D_g^{(k)}$ & $1$ & $2$ & $4$ & $8$ & $16$ & $32$\\
    \hline
    $\dim_\bQ Z_g^{(k)}$? & $1$ & $1$ & $3$ & $7$ & $15$ & $31$\\
    \hline \hline
    $\# B_h^{(k)}$ & & $1$ & $9$ & $45$ & $225$ & $1125$\\
    \hline
    $D_h^{(k)}$ & $1$ & $3$ & $8$ & $21$ & $55$ & $144$\\
    \hline
    $\dim_\bQ Z_h^{(k)}$? & $1$ & $1$ & $5$ & $15$ & $46$ & $105$\\
    \hline
   \end{tabular}
   \vspace{2mm}
 \end{table}

The above table summarizes the value of $D_f^{(k)}$ and conjectural value of $\dim_\bQ Z_f^{(k)}$ for each $f\in\{g,h\}$ and small weight $k\in\bZge{1}$ obtained by numerical experiments. The value of each iterated integral was calculated by SageMath, and linear independence was determined by the PARI/GP library of SageMath. As an observation of this table, we can check \cref{cor: dimension} numerically for $k=0,1,\dots,5$. The table also shows that there is rich relational $\bQ$-linear relations along iterated integrals on $Y_f$. The difference between $D_f^{(k)}$ and $\dim_\bQ Z_f^{(k)}$ is due to the fact that $\per\left(\left(\cH_{F_f,R_f}^{(k)}\otimes_\bQ F_f\right)^{\Gal(F_f/\bQ)}\right)\not\subseteq Z_f^{(k)}$. Also, if $f=g$, $\dim_\bQ Z_g^{(k)}$ is expected to be $2^k-1$ and the element $u\in\left(\cH_{F_f,R_f}^{(k)}\otimes_\bQ F_f\right)^{\Gal(F_f/\bQ)}$ such that $\per(u)\notin Z_g^{(k)}$ would be $u=(\logm2)^k:=\Im_\dch(0;-1;1)^k$ according to the numerical experiments.

We also have the following corollary by considering where each $\Im_f(\eta_1\cdots\eta_k)$ belongs in the direct sum decomposition in \cref{prop: direct decomposition}.

\cor \label{cor: parity}
The following conditions hold.
\begin{enumerate}[(I)]
  \item If $f=g$, we put
    \begin{align}
      Z_{g,0}^{(k)}:=\left\langle\mathrm{I}_g(\eta_1\cdots\eta_k)~\middle|~\#\left\{p~|~\eta_p=\omega_2\right\} \text{ is even}\right\rangle _\bQ,\\
      Z_{g,1}^{(k)}:=\left\langle\mathrm{I}_g(\eta_1\cdots\eta_k)~\middle|~\#\left\{p~|~\eta_p=\omega_2\right\} \text{ is odd}\right\rangle _\bQ.
    \end{align}
    Then,
    \begin{align}
      Z_{g,0}^{(k)}\subseteq\per\left(\bigoplus_{\substack{r\ge0\\k-r\text{: even}}}\Span_\bQ\left\{ b_{j_1}\cdots b_{j_r}\cdot\tau^l\otimes1~\middle|~\begin{array}{c}j_1,\dots,j_r\in\bZge{1},l\in\bZge{0},\\\sum_{i=1}^r j_i+l=k\end{array}\right\}\right),\\
      Z_{g,1}^{(k)}\subseteq\per\left(\bigoplus_{\substack{r\ge0\\k-r\text{: odd}}}\Span_\bQ\left\{ b_{j_1}\cdots b_{j_r}\cdot\tau^l\otimes\sqrt{-1}~\middle|~\begin{array}{c}j_1,\dots,j_r\in\bZge{1},l\in\bZge{0},\\\sum_{i=1}^r j_i+l=k\end{array}\right\}\right).
    \end{align}
    In particular, for each $\delta\in\{0,1\}$, it holds
    \begin{align}
      \dim_\bQ Z_{g,\delta}^{(k)}\le D_{g,\delta}^{(k)},
    \end{align}
    where
    \begin{align}
      \sum_{k=0}^\infty D_{g,0}^{(k)} t^k=\frac{1-t}{1-2t},\\
      \sum_{k=0}^\infty D_{g,1}^{(k)} t^k=\frac{t}{1-2t}.
    \end{align}

  \item If $f=h$, we put
  \begin{align}
    Z_{h,0}^{(k)}:=\left\langle\mathrm{I}_h(\eta_1\cdots\eta_k)~\middle|~\#\left\{i~|~\eta_i=\omega_4\right\} \text{ is even}\right\rangle _\bQ,\\
    Z_{h,1}^{(k)}:=\left\langle\mathrm{I}_h(\eta_1\cdots\eta_k)~\middle|~\#\left\{i~|~\eta_i=\omega_4\right\} \text{ is odd}\right\rangle _\bQ.
  \end{align}
  Then, 
  \begin{align}
    Z_{h,0}^{(k)}\subseteq\per\left(\bigoplus_{\substack{r\ge0\\k-r\text{: even}}}\Span_\bQ\left\{ b_{j_1}^{(s_1)}\cdots b_{j_r}^{(s_r)}\cdot\tau^l\otimes1~\middle|~\begin{array}{c}j_1,\dots,j_r\in\bZge{1},l\in\bZge{0},\\s_i\in\begin{cases}\{1,2\}, & j_i=1\\ \{1\}, & j_i>1\end{cases},~\sum_{i=1}^r j_i+l=k\end{array}\right\}\right),\\
    Z_{h,1}^{(k)}\subseteq\per\left(\bigoplus_{\substack{r\ge0\\k-r\text{: odd}}}\Span_\bQ\left\{ b_{j_1}^{(s_1)}\cdots b_{j_r}^{(s_r)}\cdot\tau^l\otimes\sqrt{-3}~\middle|~\begin{array}{c}j_1,\dots,j_r\in\bZge{1},l\in\bZge{0},\\s_i\in\begin{cases}\{1,2\}, & j_i=1\\ \{1\}, & j_i>1\end{cases},~\sum_{i=1}^r j_i+l=k\end{array}\right\}\right).
  \end{align}
  In particular, for each $\delta\in\{0,1\}$, it holds
  \begin{align}
    \dim_\bQ Z_{h,\delta}^{(k)}\le D_{h,\delta}^{(k)},
  \end{align}
  where
  \begin{align}
    \sum_{k=0}^\infty D_{h,0}^{(k)} t^k=\frac{1-t}{(1-3t+t^2)(1+t-t^2)},\\
    \sum_{k=0}^\infty D_{h,1}^{(k)} t^k=\frac{t(2-t)}{(1-3t+t^2)(1+t-t^2)}.
  \end{align}

\end{enumerate}

\ecor

\pf
We only show (I). Since $\Im_g(\eta_1\cdots\eta_k):=\Im_\dch(0;\theta_1\cdots\theta_k;1)$ for 
  \begin{align}
    \begin{split}
    \theta_j:=\begin{cases}
      e_0-e_{\xi_4}-e_{\xi_4^{-1}}, & \eta_j=\omega_0,\\
      -2e_1+e_{\xi_4}+e_{\xi_4^{-1}}, & \eta_j=\omega_1,\\
      \xi_4^{-1}(e_{\xi_4}-e_{\xi_4^{-1}}), & \eta_j=\omega_2,\\
      e_0, & \eta_j=\omega_3,
    \end{cases}
  \end{split}
  \end{align}
  $\Im_g(\eta_1\cdots\eta_k)$ can be written as a $\bQ$-linear sum of $u_i\otimes1$ for $u_i\in\cH_{\bQ(\sqrt{-1}),\langle2\rangle_{F_g}}^{(k)}$ if $\#\left\{p~|~\eta_p=\omega_2\right\}$ is even, and $\Im_g(\eta_1\cdots\eta_k)$ can be written as a $\bQ$-linear sum of $u_i\otimes\sqrt{-1}$ for $u_i\in\cH_{\bQ(\sqrt{-1}),\langle2\rangle_{F_g}}^{(k)}$ if $\#\left\{p~|~\eta_p=\omega_2\right\}$ is odd. The upper bound $D_{g,\delta}^{(k)}$ of the dimension $\dim_\bQ Z_{g,\delta}^{(k)}$ is given by counting the number of pairs $(j_1,\dots,j_r,l)$ satisfying each condition. (II) can be shown similarly.
\epf

The following table summarizes the value of $D_{f,\delta}^{(k)}$ and the conjectural value of $\dim_\bQ Z_{f,\delta}^{(k)}$ for each $f\in\{g,h\}$, $\delta\in\{0,1\}$, and small weight $k\in\bZge{1}$ obtained by numerical experiments. 

\begin{table}[h]
  \label{table: dimension_g}
  \centering
  \caption{dimension evaluation for $\mathrm{I}_g$}
  \begin{tabular}{|c|c|c|c|c|c|c|}
    \hline
    $k$ & $0$ & $1$ & $2$ & $3$ & $4$ & $5$\\
    \hline \hline
    $D_g^{(k)}$ & $1$ & $2$ & $4$ & $8$ & $16$ & $32$\\
    \hline
    $D_{g,0}^{(k)}$ & $1$ & $1$ & $2$ & $4$ & $8$ & $16$\\
    \hline
    $D_{g,1}^{(k)}$ & $0$ & $1$ & $2$ & $4$ & $8$ & $16$\\
    \hline \hline
    $\dim_\bQ Z_g^{(k)}$? & $1$ & $1$ & $3$ & $7$ & $15$ & $31$\\
    \hline
    $\dim_\bQ Z_{g,0}^{(k)}$? & $1$ & $0$ & $1$ & $3$ & $7$ & $15$\\
    \hline
    $\dim_\bQ Z_{g,1}^{(k)}$? & $0$ & $1$ & $2$ & $4$ & $8$ & $16$\\
    \hline
  \end{tabular}
  \vspace{2mm}
\end{table}
\begin{table}[h]
  \label{table: dimension_h}
  \centering
  \caption{dimension evaluation for $\mathrm{I}_h$}
  \begin{tabular}{|c|c|c|c|c|c|c|}
    \hline
    $k$ & $0$ & $1$ & $2$ & $3$ & $4$ & $5$\\
    \hline \hline
    $D_h^{(k)}$ & $1$ & $3$ & $8$ & $21$ & $55$ & $144$\\
    \hline
    $D_{h,0}^{(k)}$ & $1$ & $1$ & $5$ & $9$ & $30$ & $68$\\
    \hline
    $D_{h,1}^{(k)}$ & $0$ & $2$ & $3$ & $12$ & $25$ & $76$\\
    \hline \hline
    $\dim_\bQ Z_h^{(k)}$? & $1$ & $1$ & $5$ & $15$ & $46$ & $105$\\
    \hline
    $\dim_\bQ Z_{h,0}^{(k)}$? & $1$ & $0$ & $3$ & $8$ & $25$ & $53$\\
    \hline
    $\dim_\bQ Z_{h,1}^{(k)}$? & $0$ & $1$ & $2$ & $7$ & $21$ & $52$\\
    \hline
  \end{tabular}
  \vspace{2mm}
\end{table}

As an observation of these tables, we can check \cref{cor: parity} numerically for $k=0,\dots,5$. Also, it holds $\dim_\bQ Z_f^{(k)}=\dim_\bQ Z_{f,0}^{(k)}+\dim_\bQ Z_{f,1}^{(k)}$ for each $f\in\{g,h\}$, which implies that $Z_{f}^{(k)}$ may have a direct sum decomposition
\begin{align} \label{eq: direct sum decomposition for numbers}
  Z_f^{(k)}=Z_{f,0}^{(k)}\oplus Z_{f,1}^{(k)}.
\end{align}
We do not claim that the direct sum decomposition is true, but by \cref{rem: extended period conjecture}, equality (\ref{eq: direct sum decomposition for numbers}) is true if the period conjecture holds (i.e. $\per: \cH_{F_f,R_f}\rightarrow\bC$ is injective).

\rem
M. Kaneko and H. Tsumura conjectured that the $\bQ$-linear space spanned by multiple $\widetilde{T}$-values $\widetilde{T}(k_1,\dots,k_d)$ would have a direct sum decomposition by the parity of $d$ in \cite[p.5]{KaTs22}. This seems to be influenced by 
\begin{align}
  \widetilde{T}(k_1,\dots,k_d)=\mathrm{I}_g(\omega_2\omega_3^{k_1-1}\cdots\omega_2\omega_3^{k_d-1})\in \begin{cases}
    Z_{g,0}^{(k)}, & d \text{: even},\\
    Z_{g,1}^{(k)}, & d \text{: odd}.
  \end{cases}
\end{align}
in our context, and indeed the period conjecture shows their conjecture.
\erem

\subsection{Alternating multiple mixed values} \label{ss: AMMV}
The alternating multiple mixed values (AMMVs for short) are special values defined by Xu--Yan--Zhao (\cite{XuYaZh23}) and further studied by Charlton (\cite{Ch24}). AMMVs are related to the multiple $L$-values of level $4$ and include our periods in the form of $\mathrm{I}_g(\eta_1\cdots\eta_k)$. In this subsection, we observe that we can give a motivic interpretation of the AMMVs by the same method in \cref{subsec: baseext,subsec: Galoisact}.

\newcommand{\boe}{\boldsymbol{\varepsilon}}
\newcommand{\bos}{\boldsymbol{\sigma}}

\defi
For $\mathbf{k}=(k_1,\dots,k_d)\in\left(\bZge{1}\right)^d$ and $\boldsymbol{\varepsilon}=(\varepsilon_1,\dots,\varepsilon_d), \boldsymbol{\sigma}=(\sigma_1,\dots,\sigma_d)\in\{\pm1\}^d$ with $(\epsilon_d,\sigma_d)\ne(1,1)$, an alternating multiple mixed value (AMMV for short) is defined by
\begin{align}
  M_{\bos}^{\boe}(\bk):=\sum_{0<n_1<\cdots<n_d}\prod_{j=1}^d\frac{(1+\varepsilon_j(-1)^{n_j})\sigma_j^{(2m_j+1-\varepsilon_j)/4}}{n_j^{k_j}}\in\bR. \label{eq: AMMV}
\end{align}
\edefi

\rem
In \cite{XuYaZh23}, the sum of the right hand side of the \cref{eq: AMMV} was defined as $n_1>\cdots>n_d>0$, but we define as $0<n_1<\cdots<n_d$ for convenience with our setting.
\erem

As mentioned in \cite{XuYaZh23}, AMMVs have iterated integral representations as follows. Set 
\begin{align}
  \omega_0&:=\frac{\diff \lambda}{\lambda},\\
  \omega_{+1}^{-1}&:=\frac{2\diff \lambda}{1-\lambda^2}=\frac{\diff\lambda}{\lambda-1}-\frac{\diff\lambda}{\lambda+1},\\
  \omega_{-1}^{-1}&:=\frac{-2\diff \lambda}{1+\lambda^2}=\sqrt{-1}\left(\frac{\diff\lambda}{\lambda-\sqrt{-1}}-\frac{\diff\lambda}{\lambda+\sqrt{-1}}\right),\\
  \omega_{+1}^{+1}&:=\frac{2\lambda\diff\lambda}{1-\lambda^2}=\frac{\diff\lambda}{\lambda-1}+\frac{\diff\lambda}{\lambda+1},\\
   \omega_{-1}^{+1}&:=\frac{-2\lambda\diff\lambda}{1+\lambda^2}=-\frac{\diff\lambda}{\lambda-\sqrt{-1}}-\frac{\diff\lambda}{\lambda+\sqrt{-1}}
\end{align}
and
\begin{align}
  \omega^{\varepsilon_1, \varepsilon_2}_\sigma:=\begin{cases}
      -\omega_\sigma^{\varepsilon_1\varepsilon_2}, & (\sigma,\varepsilon_1,\varepsilon_2)=(-1,1,-1),\\
      \omega_\sigma^{\varepsilon_1\varepsilon_2}, & \text{otherwise}
  \end{cases}
\end{align}
for $\varepsilon_1, \varepsilon_2, \sigma\in\{\pm1\}$. Then, AMMV is written as
\begin{align} \label{eq: II rep of AMMV}
  M_{\bos}^{\boe}(\bk)=\int_\dch\omega_{\sigma_1\cdots\sigma_d}^{\varepsilon_1}\underbrace{\omega_0\cdots\omega_0}_{k_1-1}\omega_{\sigma_2\cdots\sigma_d}^{\varepsilon_2,\varepsilon_1}\underbrace{\omega_0\cdots\omega_0}_{k_2-1}\cdots\omega_{\sigma_d}^{\varepsilon_d,\varepsilon_{d-1}}\underbrace{\omega_0\cdots\omega_0}_{k_d-1}.
\end{align}

Therefore, AMMVs can be written as $\bQ(\sqrt{{-1}})$-linear sum of multiple $L$-values of level $4$. In addition, each differential form in the \cref{eq: II rep of AMMV} is $\Gal(\bQ(\sqrt{-1})/\bQ)$-invariant as same as $\theta_j$ in the \cref{eq: thetaj}. Then, we have the following proposition as an analogue of \cref{thm: 1}.

\propo
  For any $\mathbf{k}=(k_1,\dots,k_d)\in\left(\bZge{1}\right)^d$ and $\boldsymbol{\varepsilon}=(\varepsilon_1,\dots,\varepsilon_d), \boldsymbol{\sigma}=(\sigma_1,\dots,\sigma_d)\in\{\pm1\}^d$ with $(\epsilon_d,\sigma_d)\ne(1,1)$, there exists an element $M_{\bos}^{\mathfrak{m},\boe}(\bk)\in\left(\cH_{\bQ(\sqrt{-1}),\langle2\rangle_{F_g}}^{(k)}\otimes_\bQ\bQ(\sqrt{-1})\right)^{\Gal(\bQ(\sqrt{-1})/\bQ)}$ such that
  \begin{align}
    \per(M_{\bos}^{\mathfrak{m},\boe}(\bk))=M_{\bos}^{\boe}(\bk).
  \end{align}
\epropo

In particular, it holds
\begin{align}
  \dim_\bQ \mathrm{AMMV}^{(k)}\le2^k,
\end{align}
where 
\begin{align}
  \mathrm{AMMV}^{(k)}:=\mathrm{Span}_\bQ\{M_{\bos}^{\boe}(\bk)~|~k_1+\cdots+k_d=k\}.
\end{align}

\pf As mentioned above, each AMMV can be written as $\bQ(\sqrt{{-1}})$-linear sum of multiple $L$-values of level $4$. Now, each multiple $L$-value of level $4$ has a motivic interpretation as an element of $\cH_{\bQ(\sqrt{-1}),\langle2\rangle_{F_g}}$ Then, there exists an element $M_{\bos}^{\mathfrak{m},\boe}(\bk)\in\cH_{\bQ(\sqrt{-1}),\langle2\rangle_{F_g}}^{(k)}\otimes_\bQ\bQ(\sqrt{-1})$ such that
  \begin{align}
    \per(M_{\bos}^{\mathfrak{m},\boe}(\bk))=M_{\bos}^{\boe}(\bk).
  \end{align}  
  In addition, since each differential form in the \cref{eq: II rep of AMMV} is $\Gal(\bQ(\sqrt{-1})/\bQ)$-invariant, then $M_{\bos}^{\mathfrak{m},\boe}(\bk)\in\left(\cH_{\bQ(\sqrt{-1}),\langle2\rangle_{F_g}}^{(k)}\otimes_\bQ\bQ(\sqrt{-1})\right)^{\Gal(\bQ(\sqrt{-1})/\bQ)}$. In particular, since $\dim_\bQ\left(\cH_{\bQ(\sqrt{-1}),\langle2\rangle_{F_g}}^{(k)}\otimes_\bQ\bQ(\sqrt{-1})\right)^{\Gal(\bQ(\sqrt{-1})/\bQ)}=2^k$ and $$\mathrm{AMMV}^{(k)}\subseteq \per\left(\left(\cH_{\bQ(\sqrt{-1}),\langle2\rangle_{F_g}}^{(k)}\otimes_\bQ\bQ(\sqrt{-1})\right)^{\Gal(\bQ(\sqrt{-1})/\bQ)}\right),$$ we have 
\begin{align}
  \dim_\bQ \mathrm{AMMV}^{(k)}\le2^k.
\end{align}
\epf

In addition, Xu--Yan--Zhao showed 
\begin{align}
    \dim_{\bQ(\sqrt{-1})} \mathrm{AMMV}^{(k)}\otimes_\bQ\bQ(\sqrt{-1})=\dim_{\bQ(\sqrt{-1})} \mathrm{MLV}_4^{(k)}\otimes_\bQ\bQ(\sqrt{-1}),
\end{align}
where
\begin{align}
    \mathrm{MLV}_4^{(k)}:=\mathrm{Span}_\bQ\left\{L\left(\begin{matrix}
k_1,\dots,k_d\\
\alpha_1,\dots\alpha_d
\end{matrix}\right)~\middle|~k_1+\cdots+k_d=k, \alpha_1,\dots\alpha_d\in\mu_4\right\}
\end{align}
(see \cite[p.37, Theorem 6.7]{XuYaZh23}). Also, Deligne showed $$\per\left(\cH_{\bQ(\sqrt{-1}),\langle2\rangle_{F_g}}^{(k)}\right)=\mathrm{MLV}_4^{(k)}$$ (see \cite[p.37]{De10}). 
By combining these two facts, we have the following corollary.

\cor
    If the period conjecture is true (i.e. $\per: \cH_{\bQ(\sqrt{-1}),\langle2\rangle_{F_g}}\rightarrow\bC$ is injective),
\begin{align}
    \dim_\bQ \mathrm{AMMV}^{(k)}=2^k.
\end{align}
\ecor

This result is consistent with Xu--Yan--Zhao's numerical experiment (\cite[p.33, Table 1]{XuYaZh23}).

\bibliographystyle{plain}
\bibliography{bibs}

\end{document}